\documentclass[12pt]{amsart}
\usepackage{amssymb}
\usepackage{amsmath}
\usepackage{amscd}
\theoremstyle{plain}
\newtheorem{theorem}{Theorem}[section]
\newtheorem{corollary}[theorem]{Corollary}
\newtheorem{proposition}[theorem]{Proposition}
\newtheorem{lemma}[theorem]{Lemma}
{\theoremstyle{remark}

\newtheorem{remark}[theorem]{Remark}}
{\theoremstyle{definition}
\newtheorem{definition}[theorem]{Definition}
\newtheorem{example}[theorem]{Example}}

\newcommand{\benu}{\begin{enumerate}\renewcommand{\labelenumi}{{\rm (\roman{enumi})}}\renewcommand{\itemsep}{0pt}}
\newcommand{\eenu}{\end{enumerate}}

\newcommand{\N}{\mathbb{N}}
\newcommand{\Z}{\mathbb{Z}}
\newcommand{\R}{\mathbb{R}}
\newcommand{\C}{\mathbb{C}}
\newcommand{\T}{\mathbb{T}}
\newcommand{\e}{\varepsilon}
\newcommand{\cK}{{\mathcal K}}
\newcommand{\cL}{{\mathcal L}}

\newcommand{\K}{\mathbb{K}}
\newcommand{\M}{\mathbb{M}}
\newcommand{\I}{\mathbb{I}}

\newcommand{\cO}{{\mathcal O}}
\newcommand{\cT}{{\mathcal T}}
\newcommand{\cM}{{\mathcal M}}
\newcommand{\F}{{\mathcal F}}
\newcommand{\G}{{\mathcal G}}
\newcommand{\ip}[2]{\langle{#1},{#2}\rangle}
\newcommand{\bip}[2]{\big\langle{#1},{#2}\big\rangle}
\newcommand{\s}[3]{{{#1}^{#2}_{\rm{\scriptsize #3}}}}
\newcommand{\rs}[1]{{\rm{\scriptsize #1}}}
\newcommand{\case}[1]{{\medskip\noindent\underline{#1}\\}}
\DeclareMathOperator{\id}{id}
\DeclareMathOperator{\Aut}{Aut}

\DeclareMathOperator{\Ad}{Ad}
\DeclareMathOperator{\coker}{coker}
\DeclareMathOperator{\spa}{span}
\DeclareMathOperator{\supp}{supp}
\DeclareMathOperator{\cspa}{\overline{span}}

\begin{document}
\title[A class of $C^*$-algebras I, fundamental results]
{A class of \boldmath{$C^*$}-algebras generalizing 
both graph algebras and homeomorphism \boldmath{$C^*$}-algebras I, 
fundamental results}
\author[Takeshi KATSURA]{Takeshi KATSURA}
\address{
Department of Mathematical Sciences,
University of Tokyo, Komaba, Tokyo 153-8914, JAPAN}
\curraddr{
Department of Mathematics,
University of Oregon, 
Eugene, Oregon 97403-1222, U.S.A. }
\email{katsu@ms.u-tokyo.ac.jp}
\begin{abstract}
We introduce a new class of $C^*$-algebras, 
which is a generalization of both graph algebras 
and homeomorphism $C^*$-algebras. 
This class is very large and also very tractable. 
We prove the so-called gauge-invariant uniqueness theorem 
and the Cuntz-Krieger uniqueness theorem, 
and compute the $K$-groups of our algebras.
\end{abstract}

\maketitle

\setcounter{section}{-1}

\section{Introduction}

The purpose of this serial work is an introduction 
of a new class of $C^*$-algebras which contains graph algebras 
and homeomorphism $C^*$-algebras. 
Our class is very large 
so that it contains every AF-algebras \cite{Ka2}
and every Kirchberg algebras satisfying the UCT \cite{Ka4} 
as well as many simple stably projectionless $C^*$-algebras. 
At the same time, our class can be well-studied 
by using similar techniques developed in the analysis of graph algebras 
and homeomorphism $C^*$-algebras. 

Since J. Cuntz and W. Krieger introduced a class of $C^*$-algebras 
arising from finite matrices with entries $\{0,1\}$ in \cite{CK}, 
there have been many generalizations of Cuntz-Krieger algebras, 
for example, Exel-Laca algebras \cite{EL}, 
graph algebras \cite{KPRR,KPR,FLR}
and Cuntz-Pimsner algebras \cite{Pi}. 
Among others, investigation of graph algebras has rapidly progressed 
these days 
(see, for example, \cite{BPRS,BHRS,HS,DT1}), 
and many structures of graph algebras had been characterized 
in terms of graphs. 
As some authors pointed out, 
it is time to extend the techniques and results on graph algebras 
to more general $C^*$-algebras. 
Our work is one of such attempts. 
The investigation of homeomorphism $C^*$-algebras 
has also been developed 
mainly by J. Tomiyama \cite{T1, T2, T3, T4}. 
These two lines of research have several similar aspects in 
common, and our aim in this series of work is to combine and 
unify these studies in the two active fields.

In this paper, 
we associate a $C^*$-algebra 
with a quadruple $E=(E^0,E^1,d,r)$ 
where $E^0$ and $E^1$ are locally compact spaces, 
$d\colon E^1\to E^0$ is a local homeomorphism, 
and $r\colon E^1\to E^0$ is a continuous map. 
A quadruple $E=(E^0,E^1,d,r)$ is called a topological graph. 
Note that when $E^0$ is a discrete set, 
this quadruple is an ordinary (directed) graph 
and the $C^*$-algebra constructed here 
is a graph algebra of it (or its opposite graph). 
In \cite{D}, 
V. Deaconu introduced a notion of compact graphs 
and associated $C^*$-algebras with them. 
Compact graphs are particular examples of topological graphs 
and his $C^*$-algebras are isomorphic to ours. 
A triple $(E^1,d,r)$ can be considered as a generalization of 
continuous maps from $E^0$ to itself, 
and so a quadruple $E=(E^0,E^1,d,r)$ can be considered 
as a generalization of dynamical systems. 
This point of view is essential for the analysis of our $C^*$-algebras, 
and we borrow many notions from the theory of dynamical systems 
(see \cite{Ka3}). 

In the first paper of our serial work, 
we give a definition of our algebras 
and prove fundamental results on them. 
We first construct $C^*$-correspondences from topological graphs. 
This is done in Section 1 in a slightly more general form. 
Then in Section 2, we associate $C^*$-algebras 
with $C^*$-correspondences constructed from topological graphs, 
in a similar way to Cuntz-Pimsner algebras \cite{Pi}. 
We, however, point out two distinctions between our approach 
and the one for Cuntz-Pimsner algebras 
(see also the end of Section 3 in this paper). 
The first point is that left actions of our $C^*$-correspondences 
may not be injective. 
This is not allowed in \cite{Pi} 
because Cuntz-Pimsner algebras 
of $C^*$-correspondences with non-injective left actions 
often become zero 
(see \cite[Remark 1.2 (1)]{Pi}). 
Note that our algebras are obtained as relative Cuntz-Pimsner algebras 
introduced in \cite{MS}. 
The other point is that we can examine our algebras in more detail 
than Cuntz-Pimsner algebras. 
This is because our algebras are defined from ``topological'' data 
though Cuntz-Pimsner algebras are arising 
from arbitrary $C^*$-correspondences. 
In Section 3, 
we give one concrete representation of our algebras 
by using so-called Fock spaces. 
In Sections 4 and 5, we prove two fundamental theorems, namely 
the gauge-invariant uniqueness theorem (Theorem \ref{GIUT}) and 
the Cuntz-Krieger uniqueness theorem (Theorem \ref{CKUT}). 
Both of these theorems were known 
for graph algebras \cite{BPRS, RaSz, DT1, BHRS}
and for homeomorphism $C^*$-algebras \cite{T1,T2,T4} (see also \cite{ELQ}). 
We give a unified approach to these two active branches. 
In the final section, 
we prove that our algebras are always nuclear 
and satisfy the universal coefficient theorem (UCT), 
and give six-term exact sequences of $K$-groups 
which are useful to compute $K$-groups of our algebras. 
In particular, this gives a new proof 
for the computation of $K$-groups of graph algebras. 

We remark the following about our notation, 
which is different from other articles. 
We use $d$ and $r$ for denoting the ``domain'' map 
and the ``range'' map of topological graphs. 
These terms suit well when we consider quadruples 
as generalization of dynamical systems. 
However, 
our convention is the opposite one from which many papers 
on graph algebras use (for example, \cite{KPRR,KPR,FLR}). 
The author believes that 
even for graph algebras of discrete graphs, 
our convention is more natural than the one used in many papers
on graph algebras. 
One of the reasons is that 
our convention behaves well 
when one considers graphs as a kind of dynamical systems 
and graph algebras as crossed products of them. 
Another reason is that 
under our convention, 
the definition of Toeplitz families (or Cuntz-Krieger families)
$\{S_e\},\{P_v\}$ 
satisfies that 
an {\em initial} projection of $S_e$ for an edge $e$ 
is the projection $P_{d(e)}$ for the {\em domain} $d(e)$ of $e$, 
and a {\em range} projection of $S_e$ 
is a subprojection of $P_{r(e)}$ for the {\em range} $r(e)$ of $e$. 
Note that the strangeness of our definition of paths comes from 
the order of the compositions of two maps (see Section 2). 
The author is grateful to Ruy Exel 
for encouraging him to adopt this convention.

We denote the set of natural numbers 
by $\N=\{0,1,2,\ldots\}$ 
and the set of complex numbers by $\C$. 
We denote by $\T$ the group consisting of complex numbers 
whose absolute values are $1$. 
For a locally compact (Hausdorff) space $X$, 
we denote by $C(X)$ the linear space of all continuous functions on $X$. 
We define three subspaces $C_c(X)$, $C_0(X)$ and $C_b(X)$ of $C(X)$ 
by those of compactly supported functions, 
functions vanishing at infinity, and bounded functions, respectively. 

\vspace{0.5cm}

{\bf Acknowledgments.} 
The author is grateful 
to Yasuyuki Kawahigashi and Masaki Izumi 
for their encouragement, 
and to Iain Raeburn and Wojciech Szyma\'nski 
for stimulating discussions. 
He would also like to thank David Pask and Paul S. Muhly 
for valuable comments, 
and Valentin Deaconu 
for pointing out grammatical mistakes
in the earlier draft of this paper. 
This work was completed 
while the author was staying at the University of Oregon. 
He would like to thank people there, especially Huaxin Lin, 
for their warm hospitality. 
This work was partially supported by Research Fellowship 
for Young Scientists of the Japan Society for the Promotion of Science.

\section{Topological correspondences and $C^*$-correspondences}
\label{SecCorres}

In this section, 
we introduce a notion of topological correspondences and
give a method to construct $C^*$-correspondences from them. 
This method had already appeared 
in \cite{D, DKM, Sc} or other papers. 
They used completion procedure to get $C^*$-correspondences. 
However we need a concrete description of our $C^*$-correspondences. 

A (right) Hilbert $A$-module $X$ is a Banach space with
a right action of a $C^*$-algebra $A$ 
and an $A$-valued inner product $\ip{\cdot}{\cdot}$ 
satisfying that 
\benu
\item $\ip{\xi}{\eta f}=\ip{\xi}{\eta}f$, 
\item $\ip{\xi}{\eta}=\ip{\eta}{\xi}^*$, 
\item $\ip{\xi}{\xi}\geq 0$ and $\|\xi\|=\|\ip{\xi}{\xi}\|^{1/2}$, 
\eenu
for $\xi,\eta\in X$ and $f\in A$ (for the detail, see \cite{Lnc}). 
For a Hilbert $A$-module $X$, 
we denote by $\cL(X)$ the $C^*$-algebra of all adjointable operators on $X$. 
For $\xi,\eta\in X$, 
the operator $\theta_{\xi,\eta}\in\cL(X)$ is defined 
by $\theta_{\xi,\eta}(\zeta)=\xi\ip{\eta}{\zeta}$ for $\zeta\in X$.
By definition, 
$$\cK(X)=\cspa\{\theta_{\xi,\eta}\mid \xi,\eta\in X\},$$
which is an ideal of $\cL(X)$.
For $C^*$-algebras $A,B$, 
we say that $X$ is a {\em $C^*$-correspondence} from $A$ to $B$ 
when $X$ is a Hilbert $B$-module and 
a left action $\pi$ of $A$ on $X$, 
which is just a $*$-homomorphism $\pi\colon A\to \cL(X)$, 
is given. 
A $C^*$-correspondence from $A$ to itself 
is called a $C^*$-correspondence over $A$. 
When a $C^*$-correspondence $X$ from $A$ to $B$ 
and a $C^*$-correspondence $Y$ from $B$ to $C$ are given,
we can define the interior tensor product $X\otimes Y$,
which is a $C^*$-correspondence from $A$ to $C$, as follows. 
The algebraic tensor product $X\odot_B Y$ over $B$ is, by definition, 
a quotient of the algebraic tensor product $X\odot Y$ 
(as a $\C$-vector space) by the subspace 
generated by $(\xi f)\otimes\eta-\xi\otimes (\pi_Y(f)\eta)$
for $\xi\in X,\eta\in Y,f\in B$, 
where $\pi_Y\colon B\to\cL(Y)$ is the given left action.
The image of $\xi\otimes\eta\in X\odot Y$ in $X\odot_B Y$ 
is also denoted by $\xi\otimes\eta$.
We define a left action $\pi$ of $A$, 
a right action of $C$ and a $C$-valued inner product 
on $X\odot_B Y$ by 
\begin{align*}
\pi(f)(\xi\otimes\eta)&=(\pi_X(f)\xi)\otimes\eta,\\
(\xi\otimes\eta)g&=\xi\otimes(\eta g),\\
\ip{\xi\otimes\eta}{\xi'\otimes\eta'}
&=\bip{\eta}{\pi_Y(\ip{\xi}{\xi'})\eta'}
\end{align*}
for $\xi,\xi'\in X,\eta,\eta'\in Y,f\in A, g\in C$.
One can show that these operations are well-defined and 
extend to the completion of $X\odot_B Y$ 
with respect to the norm coming from the $C$-valued inner product 
defined above (see \cite[Proposition 4.5]{Lnc}).
Thus the completion of $X\odot_B Y$ is a $C^*$-correspondence 
from $A$ to $C$. 
This $C^*$-correspondence is the interior tensor product of $X$ and $Y$,
and denoted by $X\otimes Y$. 

\begin{definition}
Let $E^0$ and $E^1$ be locally compact (Hausdorff) spaces.
A map $d\colon E^1\to E^0$ is said to be {\em locally homeomorphic}
if for any $e\in E^1$, there exists a neighborhood $U$ of $e$ such that 
the restriction of $d$ to $U$ is a homeomorphism onto $d(U)$
and that $d(U)$ is a neighborhood of $d(e)$.
\end{definition}

Every local homeomorphisms are continuous and open.
If $E^0$ is discrete and 
there exists a local homeomorphism $d\colon E^1\to E^0$, 
then $E^1$ is also discrete. 

\begin{definition}\label{mvcm}
Let $E^0$ and $F^0$ be locally compact spaces. 
A {\em topological correspondence} from $E^0$ to $F^0$ 
is a triple $(E^1,d,r)$ where $E^1$ is a locally compact space, 
$d\colon E^1\to E^0$ is a local homeomorphism 
and $r\colon E^1\to F^0$ is a continuous map. 
\end{definition}

When both $d$ and $r$ are surjective local homeomorphisms, 
$(E^1,d,r)$ is called a {\em polymorphism} in \cite{AR}. 
A continuous map $\varphi\colon E^0\to F^0$ gives an example 
of topological correspondences $(E^0,\id,\varphi)$. 
More generally, a set of continuous maps 
from (possibly infinitely many) open subsets $O_i$ of $E^0$ to $F^0$ 
gives a topological correspondence 
by setting $E^1=\coprod_i O_i$ 
and defining $d$ by natural inclusions. 
Thus we consider a topological correspondence 
as a generalization of (many-valued) continuous maps. 
The pair $(E^1,d)$ defines a ``domain'' 
of a topological correspondence $(E^1,d,r)$. 
``Locally'' we can define a homeomorphism $d^{-1}$ 
from an open subset of $E^0$ to an open subset of $E^1$, and 
$r\circ d^{-1}$ gives a continuous map 
from an open subset of $E^0$ to $F^0$. 
The ``image'' of a point $v\in E^0$ under 
the topological correspondence $(E^1,d,r)$ is $r(d^{-1}(v))\subset F^0$, 
which can be more than one point (possibly infinitely many points) or 
empty. 
The ``inverse image'' of an open subset $V$ of $F^0$ is $d(r^{-1}(V))$, 
which is an open subset of $E^0$. 
When $E^0$ and $F^0$ are discrete, 
a topological correspondence $(E^1,d,r)$ 
is just a directed graph from $E^0$ to $F^0$. 

Throughout this section, 
we fix locally compact spaces $E^0,F^0$ 
and a topological correspondence $(E^1,d,r)$ 
from $E^0$ to $F^0$. 
We will see that the topological correspondence $(E^1,d,r)$ 
naturally defines a $C^*$-correspondence $C_d(E^1)$ 
from $C_0(F^0)$ to $C_0(E^0)$. 
First we define a Hilbert $C_0(E^0)$-module $C_d(E^1)$ 
by using the data of the ``domain'' $(E^1,d)$. 
For $\xi\in C(E^1)$, 
we define a map $\ip{\xi}{\xi}\colon E^0\to [0,\infty]$ by 
$\ip{\xi}{\xi}(v)=\sum_{e\in d^{-1}(v)}|\xi(e)|^2$ 
for $v\in E^0$, 
and set $\|\xi\|=\sup_{v\in E^0}\ip{\xi}{\xi}(v)^{1/2}$.
We define 
$$C_d(E^1)=\{\xi\in C(E^1)\mid \ip{\xi}{\xi}\in C_0(E^0)\}.$$
Note that we have $\|\xi\|<\infty$ for $\xi\in C_d(E^1)$. 
We will show that $C_d(E^1)$ is a Hilbert $C_0(E^0)$-module.
For $\xi,\eta\in C_d(E^1)$,
we define $\ip{\xi}{\eta}\colon E^0\to\C$ by
$$\ip{\xi}{\eta}(v)=\sum_{e\in d^{-1}(v)}\overline{\xi(e)}\eta(e)\in\C$$
for $v\in E^0$. 
This is well-defined because $\ip{\xi}{\xi}(v),\ip{\eta}{\eta}(v)<\infty$. 
To prove that $C_d(E^1)$ is a linear space, 
we need to prove $\ip{\xi}{\eta}\in C_0(E^0)$ for $\xi,\eta\in C_d(E^1)$. 
First we show that $\ip{\xi}{\eta}\in C_0(E^0)$ for $\xi,\eta\in C_c(E^1)$. 

\begin{lemma}\label{lochom1}
For any $v\in E^0$, the set $d^{-1}(v)$ has no accumulation points. 
\end{lemma}

\begin{proof}
If $d^{-1}(v)$ has an accumulation point $e\in E^1$, 
then $d|_U\colon U\to d(U)$ is not injective 
for any neighborhood $U$ of $e$.
It contradicts the fact that $d$ is locally homeomorphic.
\end{proof}

\begin{lemma}\label{lochom2}
For any compact subset $X$ of $E^1$ and any $v\in E^0$, 
there exist an open neighborhood $V$ of $v$ 
and mutually disjoint open sets $U_1,\ldots,U_n$ of $E^1$ 
such that the restriction of $d$ to $U_k$ is a homeomorphism onto $V$ 
for each $k=1,\ldots,n$ 
and that $X\cap d^{-1}(V)\subset\bigcup_{k=1}^n U_k$.
\end{lemma}

\begin{proof}
Since $X$ is compact, 
$d^{-1}(v)\cap X$ is a finite set $\{e_1,\ldots,e_n\}$ 
by Lemma \ref{lochom1}. 
Since $d$ is a local homeomorphism, 
there exist a neighborhood $V'$ of $v$ 
and a neighborhood $U'_k$ of $e_k\in E^1$ for each $k$
such that the restriction of $d$ to $U'_k$ is a homeomorphism onto $V'$.
We may assume that $U'_k$'s are mutually disjoint. 
We will show that there exists a neighborhood $V$ of $v$ 
such that $V\subset V'$ and 
$X\cap d^{-1}(V)\subset\bigcup_{k=1}^n U'_k$.
To the contrary, 
assume that for each $V\subset V'$ 
there exists $e_V\in X\cap d^{-1}(V)$ 
with $e_V\notin\bigcup_{k=1}^n U'_k$.
Since $X$ is compact, 
we can find a subnet $\{e_{V_\lambda}\}_{\lambda\in\Lambda}$ 
of the net $\{e_V\}_{V\subset V'}$ 
which converges some element $e\in X$.
We see that $d(e)=\lim_{\lambda}d(e_{V_\lambda})=\lim_{V}d(e_V)=v$.
Hence we have $e=e_k$ for some $k\in\{1,\ldots,n\}$.
Then we can find $\lambda\in\Lambda$ with $e_{V_\lambda}\in U'_k$.
This is a contradiction.
Thus we can find a neighborhood $V$ of $v$ 
such that $V\subset V'$ and 
$$X\cap d^{-1}(V)\subset \bigg(\bigcup_{k=1}^n U'_k\bigg)\cap d^{-1}(V).$$
Then $V$ and $U_k=U'_k\cap d^{-1}(V)$ are desired sets. 
\end{proof}

\begin{lemma}\label{<Cc|Cc>}
For $\xi,\eta\in C_c(E^1)$, 
we have $\ip{\xi}{\eta}\in C_c(E^0)$.
\end{lemma}

\begin{proof}
Lemma \ref{lochom1} ensures that for each $v\in E^0$, 
$\overline{\xi(e)}\eta(e)=0$ for all but finite $e\in d^{-1}(v)$.
Hence we can define $\ip{\xi}{\eta}(v)$ for $v\in E^0$, 
and we have $\ip{\xi}{\eta}\in C_c(E^0)$ 
by Lemma \ref{lochom2}.
\end{proof}

By taking $\xi=\eta$ in Lemma \ref{<Cc|Cc>}, 
we have $C_c(E^1)\subset C_d(E^1)$. 

\begin{lemma}\label{dense}
For any $\xi\in C_d(E^1)$ and any $\e>0$, 
there exists $\eta\in C_c(E^1)$ 
such that $\|\eta\|\leq\|\xi\|$ and $\|\xi-\eta\|<\e$. 
\end{lemma}

\begin{proof}
Since $\ip{\xi}{\xi}\in C_0(E^0)$, 
there exists a compact subset $Y$ of $E^0$ 
such that $\ip{\xi}{\xi}(v)<\e^2$ for $v\notin Y$.
Take $v\in Y$.
We will show that there exist a neighborhood $V_v$ of $v$ 
and a compact subset $X_v$ of $E^1$ such that 
$$\sum_{e\in d^{-1}(v')\setminus X_v}|\xi(e)|^2<\e^2,$$ 
for all $v'\in V_v$.
Since $\ip{\xi}{\xi}(v)=\sum_{e\in d^{-1}(v)}|\xi(e)|^2<\infty$,
there exist $e_1,\ldots,e_n\in d^{-1}(v)$ such that 
$$\ip{\xi}{\xi}(v)
-\sum_{k=1}^n|\xi(e_k)|^2<\frac{\e^2}{3}.$$ 
For each $k=1,\ldots,n$, 
we can find a compact neighborhood $U_k$ of $e_k$ 
such that the restriction of $d$ to $U_k$ is injective and 
$\big||\xi(e_k)|^2-|\xi(e)|^2\big|<\e^2/3n$ for $e\in U_k$.
By replacing $U_k$'s by smaller sets if necessary,
we may assume that $U_k\cap U_l=\emptyset$ for $k\neq l$. 
Since $d$ is a local homeomorphism,
$\bigcap_{k=1}^nd(U_k)$ is a neighborhood of $v$.
Hence we can find a neighborhood $V_v$ of $v$ 
such that $V_v\subset\bigcap_{k=1}^nd(U_k)$ and 
$\big|\ip{\xi}{\xi}(v')-\ip{\xi}{\xi}(v)\big|<\e^2/3$ for $v'\in V_v$.
Set $X_v=\bigcup_{k=1}^n U_k$ which is a compact subset of $E^1$.
For $v'\in V_v$, there exists a unique element $e'_k\in U_k$ 
with $d(e'_k)=v'$.
We have 
\begin{align*}
\sum_{e\in d^{-1}(v')\setminus X_v}|\xi(e)|^2
&=\ip{\xi}{\xi}(v')-\sum_{e\in d^{-1}(v')\cap X_v}|\xi(e)|^2\\
&\leq \big|\ip{\xi}{\xi}(v')-\ip{\xi}{\xi}(v)\big|
 +\bigg|\ip{\xi}{\xi}(v)-\sum_{k=1}^n|\xi(e'_k)|^2\bigg|\\
&<\frac{\e^2}{3}+\bigg|\ip{\xi}{\xi}(v)-\sum_{k=1}^n|\xi(e_k)|^2\bigg|
 +\sum_{k=1}^n\big||\xi(e_k)|^2-|\xi(e'_k)|^2\big|\\
&<\frac{\e^2}{3}+\frac{\e^2}{3}+n\frac{\e^2}{3n}=\e^2
\end{align*}
Hence we have found a neighborhood $V_v$ of $v$ 
and a compact subset $X_v$ of $E^1$ such that 
$$\sum_{e\in d^{-1}(v')\setminus X_v}|\xi(e)|^2<\e^2,$$ 
for all $v'\in V_v$.
Since $Y$ is compact, there exist $v_1,\ldots,v_m\in Y$ 
such that $Y\subset\bigcup_{k=1}^mV_{v_k}$.
Set $X=\bigcup_{k=1}^mX_{v_k}$ which is a compact subset of $E^1$. 
We can find $\eta'\in C_c(E^1)$ such that $0\leq\eta'\leq 1$ 
and $\eta'(e)=1$ for $e\in X$. 
We set $\eta=\eta'\xi\in C_c(E^1)$. 
Then we have $|\eta(e)|\leq |\xi(e)|$. 
Hence we get $\|\eta\|\leq\|\xi\|$. 
We will prove $\|\xi-\eta\|<\e$ which completes the proof.
Note that $(\xi-\eta)(e)=0$ for $e\in X$ and 
$|(\xi-\eta)(e)|\leq |\xi(e)|$ for all $e\in E^1$. 
For $v\notin Y$, we have 
$$\ip{\xi-\eta}{\xi-\eta}(v)\leq \ip{\xi}{\xi}(v)<\e^2.$$
For $v\in Y$, we can find $k$ with $v\in V_{v_k}$.
Hence we have
\begin{align*}
\ip{\xi-\eta}{\xi-\eta}(v)
&=\sum_{e\in d^{-1}(v)}|(\xi-\eta)(e)|^2
 =\sum_{e\in d^{-1}(v)\setminus X}|(\xi-\eta)(e)|^2\\
&\leq\sum_{e\in d^{-1}(v)\setminus X}|\xi(e)|^2
 \leq\sum_{e\in d^{-1}(v)\setminus X_{v_k}}|\xi(e)|^2
 <\e^2.
\end{align*}
Therefore we have $\|\xi-\eta\|<\e$. 
We are done. 
\end{proof}

By Lemma \ref{dense}, 
we see that the linear space $C_c(E^1)$ is dense in $C_d(E^1)$. 

\begin{lemma}\label{<Cd|Cd>}
For $\xi,\eta\in C_d(E^1)$, 
we have $\ip{\xi}{\eta}\in C_0(E^0)$.
\end{lemma}

\begin{proof}
For $\xi,\eta\in C_d(E^1)$, 
there exist sequences 
$\{\xi_k\}_{k\in\N}$, $\{\eta_l\}_{l\in\N}$ in $C_c(E^1)$
such that $\lim_{k\to\infty}\|\xi-\xi_k\|=0$ 
and $\lim_{l\to\infty}\|\eta-\eta_l\|=0$ 
by Lemma \ref{dense}.
By Lemma \ref{<Cc|Cc>}, $\ip{\xi_k}{\eta_l}\in C_0(E^0)$ for $k,l\in\N$.
Since 
$\|\ip{\xi_k}{\eta}-\ip{\xi_k}{\eta_l}\|\leq\|\xi_k\|\cdot\|\eta-\eta_l\|$,
we have $\ip{\xi_k}{\eta}\in C_0(E^0)$ for each $k\in\N$. 
Since 
$\|\ip{\xi}{\eta}-\ip{\xi_k}{\eta}\|\leq\|\xi-\xi_k\|\cdot\|\eta\|$,
we have $\ip{\xi}{\eta}\in C_0(E^0)$. 
\end{proof}

By Lemma \ref{<Cd|Cd>}, 
we see that $C_d(E^1)$ is a linear space. 

\begin{lemma}\label{Banach}
The linear space $C_d(E^1)$ is a Banach space
with respect to the norm $\|\cdot\|$. 
\end{lemma}

\begin{proof}
It is easy to see that $\|\cdot\|$ satisfies 
the conditions for norms. 
Take a Cauchy sequence $\{\xi_k\}_{k\in\N}$ of $C_d(E^1)$. 
Since $\sup_{e\in E}|\xi'(e)|\leq\|\xi'\|$ for $\xi'\in C_d(E^1)$, 
the sequence $\{\xi_k\}_{k\in\N}$ converges uniformly 
to some $\xi\in C(E^1)$. 
We will show that $\xi\in C_d(E^1)$ and $\{\xi_k\}_{k\in\N}$ converges 
to $\xi$ with respect to the norm $\|\cdot\|$.
For any $\e>0$, there exists $K\in\N$ such that 
$$\sum_{e\in d^{-1}(v)}|\xi_k(e)-\xi_l(e)|^2<\e^2,$$
for all $k,l\geq K$ and all $v\in E^0$.
Hence we have 
$$\sum_{e\in d^{-1}(v)}|\xi_k(e)-\xi(e)|^2\leq\e^2,$$
for all $k\geq K$ and all $v\in E^0$.
This implies that $\|\xi_k-\xi\|\leq\e$ and 
$\|\ip{\xi_k}{\xi_k}^{1/2}-\ip{\xi}{\xi}^{1/2}\|\leq\e$
for all $k\geq K$. 
Hence $\xi\in C_d(E^1)$ and $\{\xi_k\}_{k\in\N}$ converges 
to $\xi$ with respect to the norm $\|\cdot\|$. 
Thus we see that $C_d(E^1)$ is a Banach space. 
\end{proof}

Note that we have $C_c(E^1)\subset C_d(E^1)\subset C_0(E^1)$ 
and that $C_d(E^1)$ is isomorphic to the completion of $C_c(E^1)$ 
with respect to the norm $\|\cdot\|$. 

\begin{lemma}\label{rightaction}
For $\xi\in C_d(E^1)$ and $f\in C_0(E^0)$, 
we define $\xi f\colon E^1\to\C$ by $(\xi f)(e)=\xi(e)f(d(e))$.
Then $\xi f\in C_d(E^1)$ and $\ip{\eta}{\xi f}=\ip{\eta}{\xi}f$
for $\eta\in C_d(E^1)$.
\end{lemma}

\begin{proof}
It is easy to see that $\ip{\xi f}{\xi f}=\overline{f}\ip{\xi}{\xi}f$.
Hence $\xi f\in C_d(E^1)$. 
The latter part is also easy.
\end{proof}

Now we have proved the following.

\begin{proposition}\label{DefHilbmod}
The Banach space $C_d(E^1)$ is a Hilbert $C_0(E^0)$-module 
under the operations in Lemma \ref{<Cd|Cd>} and Lemma \ref{rightaction}. 
\end{proposition}

Before going further, 
we state a couple of lemmas on Hilbert modules arising from 
local homeomorphisms, 
which will be frequently used. 
Let $d$ be a local homeomorphism from $E^1$ to $E^0$, 
and $X^0$ be a closed subset of $E^0$. 
Set $X^1=d^{-1}(X^0)$ which is a closed subset of $E^1$. 
The restriction of $d$ to $X^1$ is a local homeomorphism to $X^0$.
Hence we can define a Hilbert $C_0(X^0)$-module $C_d(X^1)$ 
as in Proposition \ref{DefHilbmod}.
We use the notation $\|\cdot\|_{X}$ and $\ip{\cdot}{\cdot}_{X}$
for denoting the norm and the inner product of $C_d(X^1)$.

\begin{lemma}\label{surjection}
In the same notation as above, 
the natural map 
$C_{d}(E^1)\ni\xi\mapsto \dot{\xi}\in C_d(X^1)$ 
defined by restriction is a surjective map. 
Moreover, for $\eta\in C_d(X^1)$, 
we can find $\xi\in C_{d}(E^1)$ 
with $\dot{\xi}=\eta$ and $\|\xi\|=\|\eta\|_{X}$. 
\end{lemma}

\begin{proof}
It is easy to see that the restriction map 
$C_{d}(E^1)\ni\xi\mapsto \dot{\xi}\in C_d(X^1)$ 
is a well-defined norm-decreasing linear map. 
First we will show that for $\eta\in C_d(X^1)$,
if there exists $\xi\in C_{d}(E^1)$ with $\dot{\xi}=\eta$, 
then we can find $\zeta\in C_{d}(E^1)$ 
with $\dot{\zeta}=\eta$ and $\|\zeta\|=\|\eta\|_{X}$. 
Take $\xi\in C_{d}(E^1)$ with $\dot{\xi}=\eta$. 
Set $L=\|\eta\|_{X}$ and 
define functions $f,g\colon E^0\to [0,\infty)$ 
by $f(v)=\min\{\ip{\xi}{\xi}(v),L^2\}$
and $g(v)=L^2/\max\{L^2,\ip{\xi}{\xi}(v)\}$. 
Then we have $f\in C_0(E^0)$, $0\leq f\leq L^2$, 
$g\in C(E^0)$, $0\leq g\leq 1$ 
and $f=\ip{\xi}{\xi}g$. 
We define $\zeta\in C(E^1)$ by $\zeta(e)=\xi(e)g(d(e))^{1/2}$. 
Then we have $\ip{\zeta}{\zeta}=\ip{\xi}{\xi}g=f$.
Hence we see that $\zeta\in C_{d}(E^1)$ and $\|\zeta\|\leq L=\|\eta\|_{X}$. 
For $e\in X^1$ we have $g(d(e))=1$ 
since $\ip{\xi}{\xi}(d(e))=\ip{\eta}{\eta}_{X}(d(e))\leq L^2$. 
Hence we have $\dot{\zeta}=\dot{\xi}=\eta$ 
and this implies that $\|\zeta\|\geq\|\eta\|_{X}$. 
Thus we have shown that $\zeta\in C_{d}(E^1)$ satisfies 
that $\dot{\zeta}=\eta$ and $\|\zeta\|=\|\eta\|_{X}$. 

Next we show that the map 
$C_{d}(E^1)\ni\xi\mapsto \dot{\xi}\in C_d(X^1)$ 
is surjective. 
Take $\eta\in C_d(X^1)$ with $\eta\neq 0$. 
Set $L=\|\eta\|_{X}$. 
By Lemma \ref{dense}, 
we can find $\eta_1\in C_c(X^1)$ 
such that 
$\|\eta-\eta_1\|_{X}\leq L/2$ 
and $\|\eta_1\|_{X}\leq \|\eta\|_{X}=L$. 
Using Lemma \ref{dense} again, 
we can find $\eta_2\in C_c(X^1)$ 
such that 
$$\|(\eta-\eta_1)-\eta_2\|_{X}\leq\frac{L}{4},
\quad \mbox{ and }\quad
\|\eta_2\|\leq\|\eta-\eta_1\|_{X}\leq\frac{L}{2}.$$
Recursively, we can find $\eta_m\in C_c(X^1)$ 
such that 
$$\bigg\|\eta-\sum_{k=1}^m\eta_k\bigg\|_{X}
\leq\frac{L}{2^m}, 
\quad \mbox{ and }\quad 
\|\eta_m\|_{X}\leq\bigg\|\eta-\sum_{k=1}^{m-1}\eta_k\bigg\|_{X}
\leq\frac{L}{2^{m-1}}.$$
Then we have $\eta=\sum_{k=1}^\infty\eta_k$.
Since elements in $C_c(X^1)$ can be extend 
to elements in $C_c(E^1)$, 
we can find $\xi_k\in C_{d}(E^1)$ such that $\dot{\xi_k}=\eta_k$ 
and $\|\xi_k\|=\|\eta_k\|_{X}\leq L/2^{k-1}$ 
by the former part of this proof. 
Since $C_{d}(E^1)$ is complete, 
we can define $\xi=\sum_{k=1}^\infty\xi_k\in C_{d}(E^1)$. 
We have $\dot{\xi}=\sum_{k=1}^\infty\dot{\xi_k}=\eta$.
Hence the map 
$C_{d}(E^1)\ni\xi\mapsto \dot{\xi}\in C_d(X^1)$ 
is surjective. 
The last statement has been already proved in the argument above. 
\end{proof}

Since the restriction of $d$ to the open set $E^1\setminus X^1$ 
is a local homeomorphism to $E^0\setminus X^0$, 
we can define 
a Hilbert $C_0(E^0\setminus X^0)$-module $C_d(E^1\setminus X^1)$. 
The space $C_d(E^1\setminus X^1)$ is naturally considered 
as a subspace of $C_d(E^1)$ as 
$$C_d(E^1\setminus X^1)
=\{\xi\in C_d(E^1)\mid \xi(e)=0 \mbox{ for } e\in X^1\}
$$
Thus $C_d(E^1\setminus X^1)$ 
is a Hilbert $C_0(E^0)$-submodule of $C_d(E^1)$ and 
we see that 
$$\cK(C_d(E^1\setminus X^1))
=\cspa\{\theta_{\xi,\eta}\mid \xi,\eta\in C_{d}(E^1\setminus X^1)\}\subset 
\cK(C_d(E^1)).$$
Note that we cannot consider $\cL(C_d(E^1\setminus X^1))$
as a subspace of $\cL(C_d(E^1))$ in general. 

\begin{lemma}\label{surjection2}
For $\xi\in C_d(E^1)$, the following conditions are equivalent; 
\benu
\item $\xi\in C_d(E^1\setminus X^1)$,
\item $\dot{\xi}=0$, 
\item $\ip{\eta}{\xi}\in C_0(E^0\setminus X^0)$ for all $\eta\in C_d(E^1)$,
\item $\ip{\xi}{\xi}\in C_0(E^0\setminus X^0)$, 
\item $\xi=\xi'f$ for some $\xi'\in C_d(E^1)$ and $f\in C_0(E^0\setminus X^0)$.
\eenu
\end{lemma}

\begin{proof}
Clearly (i) is equivalent to (ii). 
(i)$\Rightarrow$(iii)$\Rightarrow$(iv) is obvious. 
For $\xi\in C_d(E^1)$ with $\ip{\xi}{\xi}\in C_0(E^0\setminus X^0)$, 
we set $f=\ip{\xi}{\xi}^{1/3}$ and 
$$\xi'(e)=\left\{
\begin{array}{ll}
\xi(e)f(d(e))^{-1} & \mbox{if } f(d(e))\neq 0\\
0 & \mbox{if } f(d(e))= 0
\end{array}\right. .$$
Then we have $f\in C_0(E^0\setminus X^0)$, 
$\xi'\in C_d(E^1)$ and $\xi=\xi'f$. 
This proves the implication (iv)$\Rightarrow$(v). 
Finally it is easy to see that (i) follows from (v).
\end{proof}

By (v) in Lemma \ref{surjection2}, 
the submodule $C_d(E^1\setminus X^1)$ of $C_d(E^1)$ 
is closed under the action of $\cL(C_d(E^1))$. 
From this fact and Lemma \ref{surjection},
we can define a $*$-homomorphism 
$\omega\colon \cL(C_d(E^1))\to\cL(C_d(X^1))$ 
by $\omega(a)\dot{\xi}=\dot{(a\xi)}$ 
for $a\in \cL(C_d(E^1))$ and $\xi\in C_d(E^1)$. 

\begin{lemma}\label{surjection3}
For $a\in \cK(C_d(E^1))$, the following conditions are equivalent; 
\benu
\item $a\in \cK(C_d(E^1\setminus X^1))$, 
\item $\omega(a)=0$, 
\item $a\xi\in C_d(E^1\setminus X^1)$ for all $\xi\in C_d(E^1)$,
\item $\ip{\eta}{a\xi}\in C_0(E^0\setminus X^0)$ 
for all $\xi,\eta\in C_d(E^1)$. 
\eenu
\end{lemma}

\begin{proof}
By the definition of $\omega$, we have (ii)$\Leftrightarrow$(iii). 
By Lemma \ref{surjection2}, we have (iii)$\Leftrightarrow$(iv). 
Clearly (i) implies (iii). 
We will prove (iii)$\Rightarrow$(i). 
Take $a\in\cK(C_d(E^1))$ 
such that $a\xi\in C_d(E^1\setminus X^1)$ for all $\xi\in C_d(E^1)$. 
There exists an approximate unit $\{u_i\}_{i\in I}$ of $\cK(C_d(E^1))$ 
such that for each $i\in I$, 
$u_i$ is a finite linear sum of elements of the form $\theta_{\xi,\eta}$. 
Since we have $a=\lim au_i$, to prove $a\in \cK(C_d(E^1\setminus X^1))$ 
it suffices to show that $a\theta_{\xi,\eta}\in \cK(C_d(E^1\setminus X^1))$ 
for arbitrary $\xi,\eta\in C_d(E^1)$. 
By the proof of Lemma \ref{surjection2}, 
we can find $\xi'\in C_d(E^1)$ and 
a positive element $f\in C_0(E^0\setminus X^0)$ such that $a\xi=\xi'f$. 
Set $\xi''=\xi'f^{1/2}$ and $\eta''=\eta f^{1/2}$. 
We have $\xi'',\eta''\in C_d(E^1\setminus X^1)$ 
and so 
$$a\theta_{\xi,\eta}=\theta_{a\xi,\eta}=\theta_{\xi'',\eta''}\in 
\cK(C_d(E^1\setminus X^1)).$$
Thus we have $a\in\cK(C_d(E^1\setminus X^1))$. 
\end{proof}

\begin{lemma}\label{surjection4}
The restriction of $\omega$ to $\cK(C_d(E^1))$ 
is a surjective map to $\cK(C_d(X^1))$, 
whose kernel is $\cK(C_d(E^1\setminus X^1))$. 
\end{lemma}

\begin{proof}
The routine computation shows that 
$\omega(\theta_{\xi,\eta})=\theta_{\dot{\xi},\dot{\eta}}$ 
for $\xi,\eta\in C_d(E^1)$. 
Hence by Lemma \ref{surjection}, 
the restriction of $\omega$ to $\cK(C_d(E^1))$ 
is a surjective map onto $\cK(C_d(X^1))$. 
We have $\cK(C_d(E^1))\cap\ker\omega=\cK(C_d(E^1\setminus X^1))$ 
by Lemma \ref{surjection3}. 
\end{proof}

There exists a $*$-homomorphism $\pi\colon C_b(E^1)\to\cL(C_d(E^1))$ 
defined by $(\pi(f)\xi)(e)=f(e)\xi(e)$ 
for $f\in C_b(E^1)$ and $\xi\in C_d(E^1)$. 
Note that $\pi$ is injective. 
We will show that for $f\in C_b(E^1)$, 
$\pi(f)\in \cK(C_d(E^1))$ if and only if $f\in C_0(E^1)$. 

\begin{lemma}\label{piK1}
If $f\in C_b(E^1)$ and $\xi_k,\eta_k\in C_d(E^1)$ $(k=1,\ldots,m)$ 
satisfy that $f=\sum_{k=1}^m\xi_k\overline{\eta_k}$
and that $\xi_k(e)\overline{\eta_k(e')}=0$
for any $k$ and any $e,e'\in E^1$ with $e\neq e'$, $d(e)=d(e')$, 
then we have $\pi(f)=\sum_{k=1}^m\theta_{\xi_k,\eta_k}$.
\end{lemma}

\begin{proof}
For $\zeta\in C_d(E^1)$ and $e\in E^1$, 
we have 
\begin{align*}
\bigg(\bigg(\sum_{k=1}^m\theta_{\xi_k,\eta_k}\bigg)\zeta\bigg)(e)
&=\sum_{k=1}^m\big(\xi_k\ip{\eta_k}{\zeta}\big)(e)
 =\sum_{k=1}^m\bigg(\xi_k(e)\sum_{e'\in d^{-1}(d(e))}
  \overline{\eta_k(e')}\zeta(e')\bigg)\\
&=\sum_{k=1}^m\sum_{e'\in d^{-1}(d(e))}
  \big(\xi_k(e)\overline{\eta_k(e')}\zeta(e')\big)
 =\sum_{k=1}^m\big(\xi_k(e)\overline{\eta_k(e)}\zeta(e)\big)\\
&=f(e)\zeta(e)
 =(\pi(f)\zeta)(e).
\end{align*}
Thus we have $\pi(f)=\sum_{k=1}^m\theta_{\xi_k,\eta_k}$.
\end{proof}

\begin{lemma}\label{piK2}
For $f\in C_c(E^1)$, 
we can find $\xi_k,\eta_k\in C_c(E^1)$ $(k=1,\ldots,m)$, 
such that $f=\sum_{k=1}^m\xi_k\overline{\eta_k}$
and that $\xi_k(e)\overline{\eta_k(e')}=0$
for any $k$ and any $e,e'\in E^1$ with $e\neq e'$ and $d(e)=d(e')$. 
\end{lemma}

\begin{proof}
We denote the support of $f$ by $X=\supp (f)$, 
which is a compact subset of $E^1$. 
Since $d$ is a local homeomorphism, 
for each $e\in X$ 
there exists an open and relatively compact neighborhood $U_e$ of $e$ 
such that the restriction of $d$ to $U_e$ is injective. 
Since $X$ is compact, 
we can find $e_1,\ldots,e_m\in X$ 
such that $X\subset\bigcup_{k=1}^mU_{e_k}$. 
Take $\zeta_1,\ldots,\zeta_m\in C_c(E^1)$ satisfying that 
$0\leq \zeta_k\leq 1$, $\supp (\zeta_k)\subset U_{e_k}$ for each $k$,
and $\sum_{k=1}^m\zeta_k(e)=1$ for all $e\in X$.
For each $k$, we define $\xi_k=f\zeta_k^{1/2}$ 
and $\eta_k=\zeta_k^{1/2}$. 
Then we have 
$$\sum_{k=1}^m\xi_k\overline{\eta_k}=f\sum_{k=1}^m\zeta_k=f.$$
For $k=1,2,\ldots,m$, 
we have $\xi_k(e)\overline{\eta_k(e')}=0$ 
for $e,e'\in E^1$ with $e\neq e'$ and $d(e)=d(e')$
because $\supp (\xi_k),\supp (\eta_k)\subset U_{e_k}$ and 
the restriction of $d$ to $U_{e_k}$ is injective.
We are done. 
\end{proof}

\begin{proposition}\label{piK}
For $f\in C_b(E^1)$, 
we have $\pi(f)\in \cK(C_d(E^1))$ if and only if $f\in C_0(E^1)$. 
\end{proposition}

\begin{proof}
For $f\in C_c(E^1)$, we have $\pi(f)\in \cK(C_d(E^1))$ 
by Lemma \ref{piK1} and Lemma \ref{piK2}. 
Hence we have $\pi(f)\in \cK(C_d(E^1))$ for every $f\in C_0(E^1)$. 
Conversely take $f\notin C_0(E^1)$. 
Then there exists $\e>0$ such that the closed set 
$$U=\{e\in E^1\mid |f(e)|\geq\e\}$$
is not compact. 
Take $\xi_1,\ldots,\xi_m,\eta_1,\ldots,\eta_m\in C_c(E^1)$
arbitrarily, 
and we will show that $\|\pi(f)-\sum_{k=1}^m\theta_{\xi_k,\eta_k}\|\geq\e$. 
Since the closed set $U$ is not compact, 
we can find $e_0\in U$ such that $e_0\notin \supp(\eta_k)$ 
for every $k=1,\ldots,m$. 
Take an open neighborhood $U_0\subset E^1$ of $e_0$ such that 
the restriction of $d$ to $U_0$ is injective 
and $U_0\cap \supp(\eta_k)=\emptyset$ for every $k=1,\ldots,m$. 
Set $\zeta\in C_c(U_0)\subset C_d(E^1)$ 
with $0\leq\zeta\leq 1$ and $\zeta(e_0)=1$. 
We have 
$$\|\zeta\|=
\sup_{v\in E^0}\bigg(\sum_{e\in d^{-1}(v)}|\zeta(e)|^2\bigg)^{1/2}
=\sup_{e\in U_0}|\zeta(e)|=1,$$
and 
\begin{align*}
\bigg\|\bigg(\pi(f)-\sum_{k=1}^m\theta_{\xi_k,\eta_k}\bigg)\zeta\bigg\|
&=\|\pi(f)\zeta\|
 =\sup_{v\in E^0}\bigg(\sum_{e\in d^{-1}(v)}|f(e)\zeta(e)|^2\bigg)^{1/2}\\
&=\sup_{e\in U}|f(e)\zeta(e)|
 \geq \e.
\end{align*}
Hence we get 
$$\bigg\|\pi(f)-\sum_{k=1}^m\theta_{\xi_k,\eta_k}\bigg\|\geq\e.$$
Since $C_c(E^1)$ is dense in $C_d(E^1)$, 
we have $\cK(C_d(E^1))=\cspa\{\theta_{\xi,\eta}\mid \xi,\eta\in C_c(E^1)\}$. 
Hence we get $\pi(f)\notin \cK(C_d(E^1))$. 
\end{proof}

\begin{remark}\label{conttrace}
We can show that $\cK(C_d(E^1))$ is a continuous trace $C^*$-algebra 
over the open subset $d(E^1)$ of $E^0$. 
For each $v\in E^0$, 
there exists a $*$-homomorphism $\cK(C_d(E^1))\to\cK(H_v)$ 
where $H_v$ is a Hilbert space whose dimension is 
the cardinality of $d^{-1}(v)$ 
(when $v\notin d(E^1)$, we set $H_v=0$). 
Hence elements of $\cK(C_d(E^1))$ 
can be considered as ``compact operator valued'' continuous 
functions on $E^0$ 
which vanish at infinity. 
Similarly elements of $\cL(C_d(E^1))$ 
can be considered as ``bounded operator valued'' 
bounded continuous functions on $E^0$. 

Each $f\in C_b(E^1)$ defines 
a topological correspondence $(E^1,d,f)$ from $E^0$ to $\C$. 
The element $f\in C_b(E^1)$ also defines $\pi(f)\in \cL(C_d(E^1))$ 
which can be considered 
as a ``diagonal operator valued'' continuous function on $E^0$. 
Thus topological correspondences from $E^0$ to $\C$
can be identified with 
``diagonal operator valued'' continuous functions on $E^0$. 
\end{remark}

So far, we only used the data of ``domain'' $(E^1,d)$ of 
the topological correspondence $(E^1,d,r)$. 
Now we will use the continuous map $r\colon E^1\to F^0$ 
to define a left action $\pi_r$ of $C_0(F^0)$ 
on the Hilbert $C_0(E^0)$-module $C_d(E^1)$. 
The continuous map $r\colon E^1\to F^0$ gives us a $*$-homomorphism 
$C_0(F^0)\ni f\to f\circ r\in C_b(E^1)$. 
Denote the composition of this map and $\pi\colon C_b(E^1)\to \cL(C_d(E^1))$ 
by $\pi_r\colon C_0(F^0)\to\cL(C_d(E^1))$. 
Explicitly, 
$(\pi_r(f)\xi)(e)=f(r(e))\xi(e)$ 
for $e\in E^1$, $f\in C_0(F^0)$ and $\xi\in C_d(E^1)$. 
In this way, 
we get a $C^*$-correspondence $C_d(E^1)$ 
from $C_0(F^0)$ to $C_0(E^0)$ 
by using a topological correspondence $(E^1,d,r)$ from $E^0$ to $F^0$. 

\begin{remark}
For $f\in C_0(F^0)$, 
we can identify 
a ``diagonal operator valued'' continuous function $\pi_r(f)$ 
with a topological correspondence 
$(E^1,d,f\circ r)$ from $E^0$ to $\C$ 
(see Remark \ref{conttrace}). 
By this observation, 
the map $\pi_r\colon C_0(F^0)\to\cL(C_d(E^1))$ 
is given by just composing 
a topological correspondence $(E^1,d,r)$. 
This observation is useful when 
we compute the $K$-groups of $\cO(E)$, 
and will be further studied in Section \ref{graphalg}
for two special examples. 
\end{remark}

\begin{lemma}\label{nondeg}
The left action $\pi_r\colon C_0(F^0)\to\cL(C_d(E^1))$ is non-degenerate. 
\end{lemma}

\begin{proof}
Take $\xi\in C_c(E^1)$, and set $K=\supp(\xi)$ 
which is a compact subset of $E^1$.
Since $r(K)$ is compact in $F^0$, 
we can find $f\in C_0(F^0)$ such that $f(v)=1$ for all $v\in r(K)$. 
Then we have $\pi_r(f)\xi=\xi$.
Since $C_c(E^1)$ is dense in $C_d(E^1)$, 
we see that 
$$\{\pi_r(f)\xi\in C_d(E^1)\mid f\in C_0(F^0),\xi\in C_d(E^1)\}$$ 
is dense in $C_d(E^1)$. 
We are done. 
\end{proof}

We define two open subsets $\s{F}{0}{sce},\s{F}{0}{fin}$ of $F^0$ by
\begin{align*}
\s{F}{0}{sce}&=\{v\in F^0\mid v\mbox{ has a neighborhood } 
V \mbox{ such that }r^{-1}(V)=\emptyset\},\\
\s{F}{0}{fin}&=\{v\in F^0\mid v\mbox{ has a neighborhood } 
V \mbox{ such that }r^{-1}(V)\mbox{ is compact}\}.
\end{align*}
We will justify the notation in Section \ref{graphalg}.
Obviously $\s{F}{0}{sce}\subset \s{F}{0}{fin}$ 
and $\s{F}{0}{sce}=F^0\setminus\overline{r(E^1)}$.
Since $\s{F}{0}{sce},\s{F}{0}{fin}$ are open, 
we can consider $C_0(\s{F}{0}{sce})$ and $C_0(\s{F}{0}{fin})$ 
as ideals of $C_0(F^0)$.

\begin{lemma}\label{O2}
Let $v\in \s{F}{0}{fin}$ and $U$ be an open subset of $E^1$ 
with $r^{-1}(v)\subset U$.
Then there exists a neighborhood $V$ of $v$ such that $r^{-1}(V)\subset U$.
\end{lemma}

\begin{proof}
Since $v\in \s{F}{0}{fin}$, 
there exists a neighborhood $V_1$ of $v$ such that $r^{-1}(V_1)$ is compact.
To derive a contradiction, suppose that 
for all neighborhood $V$ of $v$ with $V\subset V_1$, 
we can find $e_V\in r^{-1}(V)$ with $e_V\notin U$.
Since the net $\{e_V\}$ is in the compact set $r^{-1}(V_1)$, 
we can find a subnet $\{e_{V_\lambda}\}$ of $\{e_V\}$ such that 
$e_{V_\lambda}$ converges to some $e\in r^{-1}(V_1)$.
Since $U$ is open, we have $e\notin U$.
By the continuity of $r$, we have $r(e)=v$.
This contradicts the fact that $r^{-1}(v)\subset U$.
Hence we can find a neighborhood $V$ of $v$ such that $r^{-1}(V)\subset U$.
\end{proof}

\begin{lemma}\label{O2-O1}
For $v\in \s{F}{0}{fin}\setminus \s{F}{0}{sce}$ 
we have $r^{-1}(v)\neq\emptyset$.
\end{lemma}

\begin{proof}
If $v\in \s{F}{0}{fin}$ satisfies $r^{-1}(v)=\emptyset$,
then there exists a neighborhood $V$ of $v$ 
such that $r^{-1}(V)=\emptyset$ 
by Lemma \ref{O2}.
Thus we have $v\in \s{F}{0}{sce}$. 
\end{proof}

\begin{lemma}\label{r-1iscpt}
For a compact set $X\subset \s{F}{0}{fin}$, 
the subset $r^{-1}(X)$ of $E^1$ is compact.
\end{lemma}

\begin{proof}
For each $v\in X$, there exists a neighborhood $V_v$ of $v$ 
such that $r^{-1}(V_v)$ is compact.
Since $X$ is compact, we can find $v_1,\ldots,v_n\in X$ with 
$X\subset\bigcup_{i=1}^nV_{v_i}$.
Since $r^{-1}(X)\subset\bigcup_{i=1}^nr^{-1}(V_{v_i})$, 
the set $r^{-1}(X)$ is compact.
\end{proof}

\begin{proposition}\label{pir-1}
We have $\ker\pi_r=C_0(\s{F}{0}{sce})$ 
and $\pi_r^{-1}\big(\cK(C_d(E^1))\big)=C_0(\s{F}{0}{fin})$.
\end{proposition}

\begin{proof}
We have $\ker\pi_r=C_0(\s{F}{0}{sce})$ because 
\begin{align*}
f\in\ker\pi_r
&\iff f(r(e))=0 \mbox{ for all }e\in E^1\\
&\iff f(v)=0 \mbox{ for all }v\in\overline{r(E^1)}\\
&\iff f\in C_0(\s{F}{0}{sce}). 
\end{align*}
To prove the latter, 
it suffices to show that for $f\in C_0(F^0)$, 
$f\circ r\in C_0(E^1)$ if and only if $f\in C_0(\s{F}{0}{fin})$
by Proposition \ref{piK}.
If $f\in C_0(\s{F}{0}{fin})$, 
we have 
$$\{e\in E^1\mid |f(r(e))|\geq\e\}
=r^{-1}\big(\{v\in E^0\mid |f(v)|\geq\e\}\big)$$
for any $\e>0$. 
Since $\{v\in E^0\mid |f(v)|\geq\e\}$ 
is a compact subset of $\s{F}{0}{fin}$, 
Lemma \ref{r-1iscpt} shows 
that $\{e\in E^1\mid |f(r(e))|\geq\e\}$ is compact. 
Hence $f\circ r\in C_0(E^1)$. 
Now suppose $f\notin C_0(\s{F}{0}{fin})$. 
There exists $v_0\notin \s{F}{0}{fin}$ such that $|f(v_0)|>0$. 
Take $\e>0$ with $\e<|f(v_0)|$ and set $V=\{v\in E^0\mid |f(v)|\geq\e\}$. 
Then $V$ is a neighborhood of $v_0$. 
Since $v_0\notin \s{F}{0}{fin}$, 
$r^{-1}(V)$ is not compact. 
Since $\{e\in E^1\mid |f(r(e))|\geq\e\}=r^{-1}(V)$, 
we have $f\circ r\notin C_0(E^1)$. 
Therefore we have $C_0(\s{F}{0}{fin})=\pi_r^{-1}\big(\cK(C_d(E^1))\big)$.
\end{proof}

Finally we define a composition of two topological correspondences and 
prove that this relates to the interior tensor products of 
$C^*$-correspondences. 
Let $E^0,F^0,G^0$ be locally compact spaces, 
and $(E^1,d,r)$, $(F^1,d',r')$ be topological correspondences 
from $E^0$ to $F^0$ and $F^0$ to $G^0$ respectively. 
Namely, $d\colon E^1\to E^0$, $d'\colon F^1\to F^0$ 
are local homeomorphisms and 
$r\colon E^1\to F^0$, $r'\colon F^1\to G^0$ are continuous maps.
We define 
$$E^2=\{(e',e)\in F^1\times E^1\mid d'(e')=r(e)\},$$
which is a closed subset of $F^1\times E^1$.
We define a map $d''\colon E^2\to E^0$ and $r''\colon E^2\to G^0$ 
by $d''((e',e))=d(e)$ and $r''((e',e))=r'(e')$ for $(e',e)\in E^2$.

\begin{lemma}\label{compgen}
The triple $(E^2,d'',r'')$ is a topological correspondence 
from $E^0$ to $G^0$. 
\end{lemma}

\begin{proof}
Since $E^2$ is a closed subset of $F^1\times E^1$, 
it is a locally compact space.
Clearly $r''\colon E^2\to G^0$ is continuous. 
We only need to show that $d''\colon E^2\to E^0$ is locally homeomorphic. 
Take $(e',e)\in E^2$. 
There exists an open neighborhood $U'$ of $e'\in F^1$ 
such that the restriction of $d'$ to $U'$ is 
a homeomorphism onto $d'(U')$, 
and that $d'(U')$ is an open subset of $F^0$. 
We can find an open neighborhood $U$ of $e\in E^1$ 
with $U\subset r^{-1}(d'(U'))$
such that the restriction of $d$ to $U$ is 
a homeomorphism onto $d(U)$, 
and that $d(U)$ is open. 
Set $U''=E^2\cap (U'\times U)$ 
which is an open neighborhood of $(e',e)\in E^2$. 
For $v\in d(U)$, 
there exist a unique element $e_v\in U$ satisfying $d(e_v)=v$ 
and a unique element $e_{v}'\in U'$ satisfying $d'(e_{v}')=r(e_v)$. 
The map $d(U)\ni v\mapsto (e_{v}',e_v)\in U''$ 
is a continuous map which is the inverse of 
the restriction of $d''$ to $U''$.
Hence the restriction of $d''$ to $U''$ is a homeomorphism onto $d(U)$ 
which is open. 
Therefore $d''\colon E^2\to E^0$ is a local homeomorphism.
\end{proof}

The topological correspondence $(E^2,d'',r'')$ defined above 
is called the {\em composition} of two topological correspondences 
$(E^1,d,r)$ and $(F^1,d',r')$.
This composition clearly satisfies associativity. 
When $F^1=F^0$ and $d'=\id$, 
the composition of two topological correspondences 
$(E^1,d,r)$ and $(F^0,\id,r')$ is $(E^1,d,r'\circ r)$. 
We will show that the compositions of topological correspondences 
corresponds to the interior tensor products of 
$C^*$-correspondences. 
We need one lemma. 

\begin{lemma}\label{dense2}
Let $d$ be a local homeomorphism from $E^1$ to $E^0$. 
Suppose that a subset $X$ of $C_c(E^1)$ satisfies 
that for every open subset $U$ of $E^1$, every elements of $C_c(U)$ 
can be uniformly approximated by elements of $X\cap C_c(U)$. 
Then $X$ is dense in $C_d(E^1)$ with respect to the norm $\|\cdot\|$ 
\end{lemma}

\begin{proof}
Let $U$ be an open subset of $E^1$ 
such that the restriction of $d$ to $U$ is injective.
Then, 
we have $\|\xi\|=\sup_{e\in U}|\xi(e)|$ 
for $\xi\in C_c(U)\subset C_c(E^1)$. 
Hence we see that elements in $C_c(U)$ can be approximated 
by elements in $X\cap C_c(U)$ with respect to the norm $\|\cdot\|$.
By using the partition of unity, 
we can show that an arbitrary element in $C_c(E^1)$ 
is a finite sum of continuous functions 
whose supports are compact sets on which $d$ is injective. 
Hence $X$ is dense in $C_c(E^1)$ with respect to the norm $\|\cdot\|$. 
This completes the proof because $C_c(E^1)$ is dense in $C_d(E^1)$ 
with respect to the norm $\|\cdot\|$ by Lemma \ref{dense}. 
\end{proof}

\begin{proposition}\label{isom}
We have 
$$C_{d''}(E^2)\cong C_{d'}(F^1)\otimes C_d(E^1)$$
as $C^*$-correspondences from $C_0(G^0)$ to $C_0(E^0)$.
\end{proposition}

\begin{proof}
There exists a linear map 
$\varPsi\colon C_c(F^1)\odot C_c(E^1)\to C_c(E^2)$ defined by 
$$\varPsi(\xi\otimes\eta)(e',e)=\xi(e')\eta(e)\qquad 
\mbox{ for }(e',e)\in E^2.$$
Since $d'(e')=r(e)$ for $(e',e)\in E^2$, we have 
$$\varPsi((\xi f)\otimes\eta)=\varPsi(\xi\otimes (\pi_r(f)\eta))$$
for $\xi\in C_c(F^1),\eta\in C_c(E^1)$ and $f\in C_0(F^0)$.
Hence $\varPsi$ factors through 
the map $\varPsi'\colon C_c(F^1)\odot_{C_0(F^0)} C_c(E^1)\to C_c(E^2)$.
Routine computation shows that $\ip{x}{y}=\ip{\varPsi'(x)}{\varPsi'(y)}$
for $x,y\in C_c(F^1)\odot_{C_0(F^0)} C_c(E^1)$. 
Hence $\varPsi'$ extends to the isometric linear map 
$\varPsi''\colon C_{d'}(F^1)\otimes C_{d}(E^1)\to C_{d''}(E^2)$, 
which is easily shown to be a bimodule map.
To prove that $\varPsi''$ is surjective,
it suffices to show that 
the image of $\varPsi'$ is dense in $C_{d''}(E^2)$ 
because $\varPsi''$ is isometric. 
It is well-known that for each open subset $U$ of $E^2$, 
the intersection of the image of $\varPsi$ and $C_c(U)$ 
is dense in $C_c(U)\subset C_c(E^2)$ 
with respect to the sup norm. 
By Lemma \ref{dense2}, 
the image of $\varPsi'$ is dense in $C_{d''}(E^2)$ 
with respect to the norm $\|\cdot\|$.
Hence $\varPsi''$ is surjective. 
Thus $C_{d'}(F^1)\otimes C_d(E^1)$ is isomorphic to $C_{d''}(E^2)$ 
via $\varPsi''$. 
\end{proof}

\section{$C^*$-algebras arising from topological graphs}\label{graphalg}

In this section, 
we introduce a notion of topological graphs 
and give a method to define $C^*$-algebras from them. 
This construction is a generalization of ones of 
both graph algebras and homeomorphism $C^*$-algebras. 

\begin{definition}
A {\em topological graph} 
is a quadruple $E=(E^0,E^1,d,r)$ 
where $E^0,E^1$ are locally compact spaces, 
$d\colon E^1\to E^0$ is a local homeomorphism 
and $r\colon E^1\to E^0$ is a continuous map. 
\end{definition}

Note that $d,r\colon E^1\to E^0$ are not necessarily surjective nor injective.
We think that $E^0$ is a set of vertices and $E^1$ is a set of edges
and that an edge $e\in E^1$ is directed 
from its domain $d(e)\in E^0$ 
to its range $r(e)\in E^0$. 
When $E^0$ is a discrete set, then $E^1$ is also discrete. 
In this case, we call $E=(E^0,E^1,d,r)$ a {\em discrete graph}. 

For a topological graph $E=(E^0,E^1,d,r)$, 
the triple $(E^1,d,r)$ is a topological correspondence 
from $E^0$ to itself. 
Hence we can consider a quadruple $E=(E^0,E^1,d,r)$ 
as a kind of dynamical systems. 
This point of view is very important 
and we extend many notion and results 
from ordinary dynamical systems defined by homeomorphisms 
to our setting 
(see Section \ref{secCKUT} in this paper or \cite{Ka3}). 

Take a topological graph $E=(E^0,E^1,d,r)$. 
For $n=2,3,\ldots$, 
we define a space $E^n$ of {\em paths} with {\em length} $n$ by 
$$E^n=\{(e_n,\ldots,e_2,e_1)\in E^1\times \cdots\times E^1 \times E^1
\mid d(e_{k+1})=r(e_{k})\ (1\leq k\leq n-1)\}.$$ 
We define domain and range maps $d^n,r^n\colon E^n\to E^0$ by 
$d^n(e)=d(e_1)$ and $r^n(e)=r(e_n)$ for $e=(e_n,\ldots,e_1)\in E^n$. 
We write $d^1=d$ and $r^1=r$. 
We sometimes consider $E^0$ as the set of paths with length $0$. 
The domain and range maps $d^0,r^0\colon E^0\to E^0$ are defined by $d^0=r^0=\id$. 
Note that the order how to denote paths 
are the same as the one of composition of maps. 
\vspace{0.5cm}
$$(e_n,e_{n-1},\ldots,e_2,e_1)\in E^n\quad\leftrightsquigarrow\quad
\begin{CD}
\stackrel{r(e_n)}{\cdot}\hspace*{-0.3cm}
@<e_n<< \cdot
@<e_{n-1}<< \cdots\cdots @<e_{2}<< 
\cdot@<e_1<< \hspace*{-0.3cm}\stackrel{d(e_1)}{\cdot} \\
\end{CD}
\vspace*{0.5cm}$$

For $n=2,3,\ldots$, 
the triple $(E^n,d^n,r^n)$ is nothing but 
the $n$-times composition of the topological correspondence $(E^1,d,r)$. 
Hence by Lemma \ref{compgen}, 
$E^n$ is a locally compact space, 
$d^n$ is a local homeomorphism 
and $r^n$ is a continuous map for each $n\in\N$. 
From the topological correspondence $(E^1,d,r)$, 
we get a $C^*$-correspondence $C_d(E^1)$ over $C_0(E^0)$
whose left action is denoted by $\pi_r\colon C_0(E^0)\to\cL(C_d(E^1))$
as in Section \ref{SecCorres}. 
The $C^*$-correspondences $C_{d^n}(E^n)$
defined by the topological correspondence $(E^n,d^n,r^n)$ 
satisfy that for any $n,m\in\N$ 
$$C_{d^{n+m}}(E^{n+m})\cong C_{d^n}(E^n)\otimes C_{d^{m}}(E^{m})\qquad 
(\mbox{as }C^*\mbox{-correspondences over }C_0(E^0))$$ 
by Proposition \ref{isom} and for $n\geq 1$ 
$$C_{d^n}(E^n)=\cspa\{\xi_n\otimes\cdots\otimes\xi_{2}\otimes\xi_1\mid 
\xi_i\in C_{d^1}(E^1)\}.$$
Note that 
the $C^*$-correspondence $C_{d^0}(E^0)$ coincides with $C_0(E^0)$, 
and left and right actions are just multiplication. 
As long as no confusion arises, 
we omit the superscript $n$ 
and simply write $d,r$ for $d^n,r^n$. 

\begin{definition}\label{DefTpl}
Let $E=(E^0,E^1,d,r)$ be a topological graph.
A {\em Toeplitz $E$-pair} on a $C^*$-algebra $A$ 
is a pair of maps $T=(T^0,T^1)$ 
where $T^0\colon C_0(E^0)\to A$ is a $*$-homomorphism 
and $T^1\colon C_d(E^1)\to A$ is a linear map satisfying 
\benu
\item $T^1(\xi)^*T^1(\eta)=T^0(\ip{\xi}{\eta})$ for $\xi,\eta\in C_d(E^1)$, 
\item $T^0(f)T^1(\xi)=T^1(\pi_r(f)\xi)$ 
for $f\in C_0(E^0)$ and $\xi\in C_d(E^1)$.
\eenu

We denote by $\cT(E)$ the universal $C^*$-algebra
generated by a Toeplitz $E$-pair.
\end{definition}

For a Toeplitz $E$-pair $T=(T^0,T^1)$, 
we see that $\|T^0(f)\|\leq\|f\|$ and $\|T^1(\xi)\|\leq\|\xi\|$ 
because $T^0$ is a $*$-homomorphism and 
$$\|T^1(\xi)\|^2=\|T^1(\xi)^*T^1(\xi)\|=\|T^0(\ip{\xi}{\xi})\|
\leq\|\ip{\xi}{\xi}\|=\|\xi\|^2.$$
Hence the universal $C^*$-algebra $\cT(E)$ 
generated by a Toeplitz $E$-pair exists
(see Section \ref{Fock} for a concrete construction of $\cT(E)$). 
We have $T^1(\xi)T^0(f)=T^1(\xi f)$ for $f\in C_0(E^0)$ and $\xi\in C_d(E^1)$
because $(T^1(\xi)T^0(f)-T^1(\xi f))^*(T^1(\xi)T^0(f)-T^1(\xi f))=0$ 
by the condition (i) above. 
We write $C^*(T)$ for denoting the $C^*$-algebra 
generated by the images of the maps $T^0$ and $T^1$. 
Let $n$ be an integer greater than $1$, 
and $\xi_1,\ldots,\xi_n$, $\eta_1,\ldots,\eta_n$ 
be elements of $C_d(E^1)$. 
Set $\xi=\xi_n\otimes\cdots\otimes\xi_1,\ 
\eta=\eta_n\otimes\cdots\otimes\eta_1\in C_d(E^n)$. 
By using condition (i) and (ii) in Definition \ref{DefTpl}, 
we can prove 
$$\big(T^1(\xi_n)\cdots T^1(\xi_1)\big)^*
\big(T^1(\eta_n)\cdots T^1(\eta_1)\big)=T^0(\ip{\xi}{\eta}).$$
Therefore 
we can define a norm-decreasing linear map $T^n\colon C_d(E^n)\to C^*(T)$ 
by $T^n(\xi)=T^1(\xi_n)\cdots T^1(\xi_1)$
for $\xi=\xi_n\otimes\cdots\otimes\xi_1\in C_d(E^n)$. 
For $n\in\N$, 
we define a linear map $\varPhi^n$ from 
$\spa\{\theta_{\xi,\eta}\in\cK(C_d(E^n))\mid \xi,\eta\in C_d(E^n)\}$
to $C^*(T)$ by 
$$\varPhi^n(\theta_{\xi,\eta})=T^n(\xi)T^n(\eta)^*.$$
One can check that this map is a well-defined norm-decreasing $*$-homomorphism 
(see \cite[Lemma 3.2]{Pi} or \cite[Lemma 2.1]{KPW}).
Hence it uniquely extends to a $*$-homomorphism 
$\varPhi^n\colon \cK(C_d(E^n))\to C^*(T)$. 
Note that $\varPhi^0=T^0$ if we identify $\cK(C_d(E^0))$ 
with $C_0(E^0)$ in the natural way.
We summarize properties of $T^n$ and $\varPhi^n$ in the following lemma.
The proof is left to the reader.

\begin{lemma}\label{Tn}
Let $E=(E^0,E^1,d,r)$ be a topological graph 
and $T=(T^0,T^1)$ be a Toeplitz $E$-pair.
Then the maps $T^n\colon C_d(E^n)\to C^*(T)$ 
and $\varPhi^n\colon \cK(C_d(E^n))\to C^*(T)$ 
defined above satisfy the following 
$(n,m\in\N, \xi,\zeta\in C_d(E^n),\eta\in C_d(E^m), 
f\in C_0(E^0), x\in\cK(C_d(E^n)))$: 
\benu
\item $T^n(\xi)T^m(\eta)=T^{n+m}(\xi\otimes\eta)$, 
\item $T^n(\zeta)^*T^n(\xi)=T^0(\ip{\zeta}{\xi})$, 
\item $T^0(f)T^n(\xi)=T^n(\pi_{r^n}(f)\xi)$, 
\item $T^0(f)\varPhi^n(x)=\varPhi^n(\pi_{r^n}(f)x)$, 
\item $\varPhi^n(x)T^n(\xi)=T^n(x\xi)$. 
\eenu
\end{lemma}

We say that a Toeplitz $E$-pair $T=(T^0,T^1)$ is {\em injective} 
if $T^0$ is injective. 
It is easy to see that for an injective Toeplitz $E$-pair $T$, 
$T^n$ and $\varPhi^n$ are isometric for all $n\in\N$. 

\begin{lemma}\label{k=m-n}
Let $n,m\in \N$ be integers with $n<m$, 
and $\xi\in C_d(E^n)$, $\eta\in C_d(E^m)$. 
Then we have $T^{k}(\zeta)=T^n(\xi)^*T^m(\eta)$ 
where $k=m-n\in \N$ and 
$\zeta\in C_d(E^{k})$ is defined by 
$$\zeta(e)
=\sum_{\stackrel
{\mbox{$\scriptstyle e'\in E^n$}}
{\mbox{$\scriptstyle d(e')=r(e)$}}}
\overline{\xi(e')}\eta(e',e)\qquad (e\in E^k).$$
\end{lemma}

\begin{proof}
Take $\xi,\eta_1\in C_d(E^n)$ and $\eta_2\in C_d(E^{k})$, 
and set $\eta=\eta_1\otimes\eta_2\in C_d(E^m)$. 
The element $\zeta\in C_d(E^{k})$ defined by the above equation 
satisfies $\zeta=\pi_r(\ip{\xi}{\eta_1})\eta_2$ 
because we have 
$$\zeta(e)
=\sum_{\stackrel{\mbox{$\scriptstyle e'\in E^n$}}
{\mbox{$\scriptstyle d(e')=r(e)$}}}
\overline{\xi(e')}\eta_1(e')\eta_2(e)
=\ip{\xi}{\eta_1}(r(e))\eta_2(e)
=\big(\pi_r(\ip{\xi}{\eta_1})\eta_2\big)(e),$$
for $e\in E^k$.
By Lemma \ref{Tn}, 
we get 
\begin{align*}
T^n(\xi)^*T^{m}(\eta)
&=T^n(\xi)^*T^{n}(\eta_1)T^{k}(\eta_1)\\
&=T^0(\ip{\xi}{\eta_1})T^{k}(\eta_1)\\
&=T^k(\pi_r(\ip{\xi}{\eta_1})\eta_2)\\
&=T^k(\zeta).
\end{align*}
The set of linear combinations 
of elements of the form $\eta_1\otimes\eta_2$ 
($\eta_1\in C_d(E^n)$, $\eta_2\in C_d(E^{k})$) is dense in $C_d(E^m)$. 
Hence the equation holds for all $\xi\in C_d(E^n)$ and 
$\eta\in C_d(E^m)$. 
\end{proof}

By the above lemma, we have 
$$C^*(T)=\cspa\{T^n(\xi)T^m(\eta)^*\mid 
\xi\in C_d(E^n),\ \eta\in C_d(E^m),\ n,m\in\N\}.$$
Combining this fact with Lemma \ref{nondeg}, 
we can easily show that the hereditary $C^*$-algebra generated 
by $T^0(C_0(E^0))\subset C^*(T)$ is $C^*(T)$. 
From this fact, we get the following. 

\begin{proposition}\label{strict}
A net $\{u_i\}$ 
in the multiplier algebra $\cM(C^*(T))$ of $C^*(T)$ 
converges in the strict topology 
if and only if $u_iT_0(f)$ and $T_0(f)u_i$ 
converge to elements in $C^*(T)$ 
in the norm topology for every $f\in C_0(E^0)$. 
\end{proposition}

To introduce Cuntz-Krieger $E$-pairs of a graph $E$, 
we need the following notion. 

\begin{definition}
Let $E=(E^0,E^1,d,r)$ be a topological graph.
We define three open subsets 
$\s{E}{0}{sce},\s{E}{0}{fin}$ and $\s{E}{0}{rg}$ of $E^0$ by
$\s{E}{0}{sce}=E^0\setminus\overline{r(E^1)}$,
\begin{align*}
\s{E}{0}{fin}=\{v\in E^0\mid\mbox{ there exists}&\mbox{ a neighborhood } 
V \mbox{ of } v\\
&\mbox{ such that }
r^{-1}(V)\subset E^1 \mbox{ is compact}\},
\end{align*}
and $\s{E}{0}{rg}=\s{E}{0}{fin}\setminus\overline{\s{E}{0}{sce}}$.
We define two closed subsets $\s{E}{0}{inf}$ and $\s{E}{0}{sg}$ of $E^0$ by 
$\s{E}{0}{inf}=E^0\setminus \s{E}{0}{fin}$ and 
$\s{E}{0}{sg}=E^0\setminus \s{E}{0}{rg}$. 
\end{definition}

A vertex in $\s{E}{0}{sce}$ is called a {\em source}.
When $E$ is a discrete graph, 
$\s{E}{0}{fin}$ is the set of vertices 
which receive finitely many edges, 
while $\s{E}{0}{inf}$ is the set of vertices 
which receive infinitely many edges. 
A vertex in $\s{E}{0}{rg}$ is said to be {\em regular}, 
and a vertex in $\s{E}{0}{sg}$ is said to be {\em singular}. 
We see that $\s{E}{0}{sg}=\overline{\s{E}{0}{sce}}\cup\s{E}{0}{inf}$.
Since $\s{E}{0}{sce}\subset \s{E}{0}{fin}$, 
we have $\s{E}{0}{sce}\cap\s{E}{0}{inf}=\emptyset$. 
However it may happen that 
$\overline{\s{E}{0}{sce}}\cap\s{E}{0}{inf}\neq\emptyset$, 
as the following example shows. 

\begin{example}
Define a topological graph $E=(E^0,E^1,d,r)$ 
by $E^0=\R$, $E^1=(0,\infty)\subset\R$ 
and $d,r$ are natural embeddings. 
Then we have $\s{E}{0}{sce}=(-\infty,0), \s{E}{0}{fin}=\R\setminus\{0\}$. 
Hence $\s{E}{0}{inf}=\{0\}$ and $\overline{\s{E}{0}{sce}}=(-\infty,0]$ 
have a non-empty intersection. 
The set of regular vertices is $\s{E}{0}{rg}=(0,\infty)$. 
\end{example}

\begin{proposition}\label{E0r}
For $v\in E^0$, we have $v\in\s{E}{0}{rg}$ 
if and only if there exists a neighborhood $V$ of $v$ such that 
$r^{-1}(V)$ is compact and $r(r^{-1}(V))=V$.
\end{proposition}

\begin{proof}
If there exists a neighborhood $V$ of $v$ such that 
$r^{-1}(V)$ is compact and $r(r^{-1}(V))=V$, 
then $v\in\s{E}{0}{fin}$ and $V\cap \s{E}{0}{sce}=\emptyset$.
Hence we get $v\in\s{E}{0}{rg}$. 
Conversely if $v\in\s{E}{0}{rg}$, 
then a compact neighborhood $V$ of $v$ with $V\subset\s{E}{0}{rg}$ 
satisfies that $r^{-1}(V)$ is compact by Lemma \ref{r-1iscpt}
and that $r(r^{-1}(V))=V$ by Lemma \ref{O2-O1}. 
\end{proof}

Proposition \ref{E0r} means that the open set $\s{E}{0}{rg}$ 
is largest among open subsets $U$ of $E^0$ 
such that the restriction of $r$ to $r^{-1}(U)$ 
is a proper surjection onto $U$. 
Note that for $v\in\s{E}{0}{rg}$ 
$r^{-1}(v)$ is a non-empty compact set 
by Proposition \ref{E0r}. 

\begin{definition}
Let $E=(E^0,E^1,d,r)$ be a topological graph.
A Toeplitz $E$-pair $T=(T^0,T^1)$ is called a {\em Cuntz-Krieger $E$-pair} 
if $T^0(f)=\varPhi^1(\pi_r(f))$ holds for all $f\in C_0(\s{E}{0}{rg})$.
\end{definition}

Note that 
the restriction of $\pi_r$ to $C_0(\s{E}{0}{rg})$ is 
an injection into $\cK(C_d(E^1))$ 
by Proposition \ref{pir-1}. 

\begin{definition}
We denote by $\cO(E)$ the universal $C^*$-algebra
generated by a Cuntz-Krieger $E$-pair $t=(t^0,t^1)$.
\end{definition}

For $n\in\N$, 
we write $t^n\colon C_d(E^n)\to \cO(E)$ 
and $\varphi^n\colon \cK(C_d(E^n))\to \cO(E)$ 
for denoting the maps corresponding to $T^n$ and $\varPhi^n$.
For a Cuntz-Krieger $E$-pair $T=(T^0,T^1)$, 
we denote by $\rho_{T}$ 
the unique surjection from $\cO(E)$ to $C^*(T)$ 
satisfying $\rho_{T}\circ t^i=T^i$ for $i=0,1$.
The map $\rho_{T}$ satisfies $\rho_{T}\circ t^n=T^n$ and 
$\rho_{T}\circ \varphi^n=\varPhi^n$ 
for all $n\in\N$.
At this point, 
we do not know whether there exists 
an injective Cuntz-Krieger $E$-pair, 
or even an injective Toeplitz $E$-pair. 
In the next section, 
we will construct 
one injective Cuntz-Krieger $E$-pair $\tau=(\tau^0,\tau^1)$.  
This implies that the universal Cuntz-Krieger $E$-pair $t=(t^0,t^1)$ 
is injective. 
The next lemma may help us to understand 
a role of $\s{E}{0}{rg}$ 
and the definition of Cuntz-Krieger $E$-pairs.

\begin{proposition}\label{T=Pp}
Let $T$ be an injective Toeplitz $E$-pair. 
If $f\in C_0(E^0)$ satisfies $T^0(f)\in\varPhi^1\big(\cK(C_d(E^1))\big)$, 
then $f\in C_0(\s{E}{0}{rg})$ and $T^0(f)=\varPhi^1(\pi_r(f))$. 
\end{proposition}

\begin{proof}
Let $f$ be an element of $C_0(E^0)$ 
satisfying $T^0(f)\in\varPhi^1\big(\cK(C_d(E^1))\big)$. 
Take $x\in\cK(C_d(E^1))$ with $T^0(f)=\varPhi^1(x)$.
For $\xi\in C_d(E^1)$, we see that 
$$T^1(\pi_r(f)\xi)=T^0(f)T^1(\xi)
=\varPhi^1(x)T^1(\xi)=T^1(x\xi).$$
Since $T^1$ is injective, 
we get $\pi_r(f)=x$. 
Thus we have $T^0(f)=\varPhi^1(\pi_r(f))$. 
By Proposition \ref{pir-1}, we have $f\in C_0(\s{E}{0}{fin})$. 
We will show that $f\in C_0(\s{E}{0}{rg})$. 
To derive a contradiction, 
assume that there exists $v\notin \s{E}{0}{rg}$ such that $f(v)\neq 0$. 
Since $f\in C_0(\s{E}{0}{fin})$, 
we see that $v\in\overline{\s{E}{0}{sce}}$. 
Hence we can find $v'\in \s{E}{0}{sce}$ such that $f(v')\neq 0$.
The element $v'\in \s{E}{0}{sce}$ has a neighborhood $V$ 
with $r^{-1}(V)=\emptyset$. 
Take $g\in C_0(V)$ with $g(v')\neq 0$.
Then we have $fg\neq 0$ and 
$$T^0(fg)=T^0(f)T^0(g)=\varPhi^1(\pi_r(f))T^0(g)
=\varPhi^1(\pi_r(f)\pi_r(g))=0,$$
because $\pi_r(g)=0$. 
This contradicts the fact that $T^0$ is injective.
Therefore $f\in C_0(\s{E}{0}{rg})$. 
\end{proof}

\begin{remark}
In our sequel \cite{Ka2},
we will show that $\cO(E)$ 
is the smallest $C^*$-algebra which is generated 
by an {\em injective} Toeplitz $E$-pair which admits a gauge action
(which means that there exists an automorphism $\beta'_z$ of $C^*(T)$ 
with $\beta'_z(T^0(f))=T^0(f)$ and $\beta'_z(T^1(\xi))=zT^1(\xi)$
for every $z\in\T$).
\end{remark}

We give two fundamental examples of topological graphs 
and $C^*$-algebras associated with them. 
More elaborated examples can be found in \cite{Ka2}.

\vspace{0.3cm}
\noindent\underline{\bf Example 1 (graph algebras)}

When $E^0$ is discrete, $E^1$ is also discrete 
and $E=(E^0,E^1,d,r)$ becomes an ordinary (directed) graph. 
We have 
\begin{align*}
\s{E}{0}{fin}&=\{v\in E^0\mid r^{-1}(v)\mbox{ is a finite set}\},\\ 
\s{E}{0}{sce}&=\{v\in E^0\mid r^{-1}(v)=\emptyset\},\\ 
\s{E}{0}{rg}&=\{v\in E^0\mid r^{-1}(v)
\mbox{ is a non-empty finite set}\}. 
\end{align*}
For a Toeplitz $E$-pair $(T^0,T^1)$, 
define $P_v=T^0(\delta_v)$ for $v\in E^0$ 
where $\delta_v\in C_0(E^0)$ is a characteristic function on $\{v\}$. 
Then $\{P_v\}_{v\in E^0}$ is a family of mutually orthogonal 
projections. 
Similarly set $S_e=T^1(\delta_e)$ for $e\in E^1$. 
Then the family $(\{P_v\}_{v\in E^0},\{S_e\}_{e\in E^1})$ is 
a Toeplitz-Cuntz-Krieger $\bar{E}$-family in the sense of \cite{FR}
where $\bar{E}$ is the opposite graph of $E$, 
that is, the set of vertices and edges of $\bar{E}$ 
are the same as those of $E$, but the range map of $\bar{E}$ is $d$ 
and the source map of $\bar{E}$ is $r$. 
Conversely, from a Toeplitz-Cuntz-Krieger $\bar{E}$-family 
$(\{P_v\}_{v\in E^0},\{S_e\}_{e\in E^1})$, 
we can define a Toeplitz $E$-pair $(T^0,T^1)$ by 
$T^0(f)=\sum_{v\in E^0}f(v)P_v$ and $T^1(\xi)=\sum_{e\in E^1}\xi(e)S_e$. 
Thus there exists a one-to-one correspondence 
between the set of Toeplitz $E$-pairs 
and the set of Toeplitz-Cuntz-Krieger $\bar{E}$-families. 
Under this correspondence, 
Cuntz-Krieger $E$-pairs correspond exactly to 
Cuntz-Krieger $\bar{E}$-families in the sense of \cite{FLR}. 
Thus $\cO(E)$ is isomorphic to the graph algebra of the graph $\bar{E}$. 

We can describe $\cK(C_d(E^1))$ explicitly in this case. 
For $e,e'\in E^1$, 
we define $u_{e,e'}=\theta_{\delta_e,\delta_{e'}}\in\cK(C_d(E^1))$. 
Then we have 
$$\cK(C_d(E^1))=\cspa\{u_{e,e'}\mid e,e'\in E^1\}.$$

\begin{lemma}\label{disgra1}
We have $u_{e,e'}\neq 0$ if and only if $d(e)=d(e')$. 
\end{lemma}

\begin{proof}
When $e''\neq e'$, 
we have $u_{e,e'}(\delta_{e''})=0$  
because $\ip{\delta_{e'}}{\delta_{e''}}=0$. 
We have $u_{e,e'}(\delta_{e'})=\delta_e\delta_{d(e')}$, 
and $\delta_e\delta_{d(e')}$ is non-zero if and only if $d(e)=d(e')$ 
(and in this case $\delta_e\delta_{d(e')}=\delta_e$). 
Hence $u_{e,e'}\neq 0$ if and only if $d(e)=d(e')$. 
\end{proof}

\begin{lemma}\label{disgra2}
For $e_1,e_1',e_2,e_2'\in E^1$ with $d(e_1)=d(e_1'),d(e_2)=d(e_2')$, 
we have 
$$u_{e_1,e_1'}u_{e_2,e_2'}=
\left\{\begin{array}{ll}
u_{e_1,e_2'}& ({\rm if}\ e_1'=e_2)\\
0& ({\rm if}\ e_1'\neq e_2)\\
\end{array}\right.$$
\end{lemma}

\begin{proof}
Clear by the computation in the proof of Lemma \ref{disgra1}. 
\end{proof}

For $v\in E^0$, 
define $K_v=\cspa\{u_{e,e'}\mid e,e'\in d^{-1}(v)\}$. 
By Lemma \ref{disgra1} and Lemma \ref{disgra2}, 
we have the following.

\begin{lemma}\label{disgra3}
\benu
\item If $d^{-1}(v)=\emptyset$, then $K_v=0$. 
\item If $d^{-1}(v)$ is infinite, then $K_v\cong\K$. 
\item If $d^{-1}(v)$ consists of $n$ edges, 
then $K_v\cong\M_{n}$.
\item For distinct $v,v'\in E^0$, 
$K_v$ and $K_{v'}$ are orthogonal to each others. 
\item $\cK(C_d(E^1))=\bigoplus_{v\in E^0}K_v$. 
\eenu
\end{lemma}

By Lemma \ref{disgra3}, 
it is easy to see that 
there exists a natural isomorphism between 
the $K$-groups of $\cK(C_d(E^1))$ and the ones of $C_0(d(E^1))$. 
This is the case for general topological graphs $E$ 
because the Hilbert module $C_d(E^1)$ gives 
a strong Morita equivalence between $\cK(C_d(E^1))$ 
and $C_0(d(E^1))\subset C_0(E^0)$ (see \cite{E}). 
The map $\pi_r\colon C_0(\s{E}{0}{fin})\to \cK(C_d(E^1))$ 
can be described as 
$$\pi_r(f)=\sum_{e\in E^1}f(r(e))u_{e,e}
=\bigoplus_{v\in E^0}\hat{r}(f)(v),$$
where $\hat{r}(f)(v)=\sum_{e\in d^{-1}(v)}f(r(e))u_{e,e}$ 
is a diagonal operator of $K_v$ for $v\in E^0$. 
This formula makes sense for the map 
$\pi_r\colon C_0(E^0)\to \cL(C_d(E^1))$, 
where infinite sums converge in strict topology. 

\vspace{0.3cm}
\noindent\underline{\bf Example 2 (homeomorphism \boldmath{$C^*$}-algebras)}

Take a topological dynamical system $\Sigma=(X,\sigma)$ 
where $X$ is a compact space 
and $\sigma\colon X\to X$ is a homeomorphism. 
We can define an automorphism $\alpha$ of $C(X)$ 
by $\alpha(f)(x)=f(\sigma^{-1}x)$. 
The crossed product $C(X)\rtimes_{\alpha}\Z$ 
is called a {\em homeomorphism $C^*$-algebra} 
and denoted by $A(\Sigma)$ in \cite{T3,T4}. 
A representation of $A(\Sigma)$ corresponds bijectively to 
a covariant representation $\{\pi,u\}$ 
of the topological dynamical system $\Sigma$ 
where $\pi$ is a representation of $C(X)$ 
and $u$ is a unitary satisfying that $\pi(\alpha(f))=u\pi(f)u^*$. 
A homeomorphism $\sigma$ defines a topological correspondence 
on $X$, hence we get a topological graph $E_\Sigma$ from $\Sigma$. 
We will see that the homeomorphism $C^*$-algebra $A(\Sigma)$ is 
isomorphic to $\cO(E_\Sigma)$. 

We treat a more general setting, 
namely when $X$ is a locally compact space, 
and $\sigma$ is a proper continuous map from $X$ to $X$. 
Define a topological graph $E_\Sigma=(E_\Sigma^0,E_\Sigma^1,d,r)$ by 
$E_\Sigma^0=E_\Sigma^1=X$, $d=\id_X$ and $r=\sigma$. 
For a natural number $n\geq 2$, 
$$E_\Sigma^n=\{(x_n,\ldots,x_1)\in X\times\cdots\times X \mid 
x_k=\sigma(x_{k-1}) \mbox{ for } k=2,\ldots,n\}$$
is isomorphic to $X$ by $E_\Sigma^n\ni (x_n,\ldots,x_1)\mapsto x_1\in X$. 
We will identify $E_\Sigma^n$ with $X$ by this map. 
Under this identification, 
we see that $d^n=\id_X$ and $r^n=\sigma^n$. 
In other words, 
$(X,\id_X,\sigma^n)$ is the $n$-times composition of 
the topological correspondence $(X,\id_X,\sigma)$. 
We identify $C_{d^n}(E_\Sigma^n)$ with $C_0(X)$ for every $n\in\N$
as (right) Hilbert $C_0(X)$-modules. 
We define an endomorphism $\hat{\sigma}\colon C_0(X)\to C_0(X)$ 
by $\hat{\sigma}(f)=f\circ\sigma$. 
If we identify $\cK(C_d(E_\Sigma^1))$ with $C_0(X)$, 
then the map $\hat{\sigma}$ coincides 
with the left action 
$\pi_r\colon C_0(X)\to\cK(C_d(E_\Sigma^1))\subset \cL(C_d(E_\Sigma^1))$ 
defined by $r\ (=\sigma)$. 
Let us take a Toeplitz $E_\Sigma$-pair $T=(T^0,T^1)$. 

\begin{lemma}\label{tga1}
For $\xi,\eta\in C_0(X)$ and $n,m\in\N$, 
we have the following. 
\benu 
\item $T^n(\xi)^*T^n(\eta)=T^0(\overline{\xi}\eta)$.
\item $T^n(\xi)T^m(\eta)=T^{n+m}(\hat{\sigma}^m(\xi)\eta)$. 
\item $T^n(\xi)^*T^m(\eta)=T^{m-n}(\hat{\sigma}^{m-n}(\overline{\xi})\eta)$
when $n\leq m$. 
\eenu
\end{lemma}

\begin{proof}
\benu 
\item Clear by $d^n=\id_X$.
\item We have $T^n(\xi)T^m(\eta)=T^{n+m}(\xi\otimes\eta)$ 
and 
$$\xi\otimes\eta(x)=
\xi\otimes\eta(\sigma^{n+m-1}(x),\ldots,\sigma(x),x)
=\xi(\sigma^{m}(x))\eta(x)=\big(\hat{\sigma}^m(\xi)\eta\big)(x).$$
Hence we have $T^n(\xi)T^m(\eta)=T^{n+m}(\hat{\sigma}^m(\xi)\eta)$. 
\item Since the set 
$\{\hat{\sigma}^{m-n}(\eta_1)\eta_2\mid \eta_1,\eta_2\in C_0(X)\}$
is dense in $C_0(X)$, 
it suffices to show the equation for 
$\eta=\hat{\sigma}^{m-n}(\eta_1)\eta_2$ where $\eta_1,\eta_2\in C_0(X)$. 
We have $T^m(\eta)=T^n(\eta_1)T^{m-n}(\eta_2)$ by (ii).
Hence 
\begin{align*}
T^n(\xi)^*T^m(\eta)
&=T^n(\xi)^*T^n(\eta_1)T^{m-n}(\eta_2)
=T^0(\overline{\xi}\eta_1)T^{m-n}(\eta_2)\\
&=T^{m-n}(\hat{\sigma}^{m-n}(\overline{\xi}\eta_1)\eta_2)
=T^{m-n}(\hat{\sigma}^{m-n}(\overline{\xi})\eta).
\end{align*}
\eenu
\end{proof}

\begin{proposition}\label{tga2}
Take an approximate unit $\{u_i\}$ of $C_0(X)$. 
Then the net $\{T^1(u_i)\}$ in $C^*(T)$ converges strictly 
to an element $U_T\in\cM(C^*(T))$ 
satisfying $U_TT^0(f)=T^1(f)$ and $T^0(f)U_T=T^1(\hat{\sigma}(f))$ 
for $f\in C_0(X)$. 
\end{proposition}

\begin{proof}
For $f\in C_0(X)$, we have $T^1(u_i)T^0(f)=T^1(u_i f)$ 
which converges to $T^1(f)$ in the norm topology. 
We also have that $T^0(f)T^1(u_i)=T^1(\hat{\sigma}(f)u_i)$ 
converges to $T^1(\hat{\sigma}(f))$ in the norm topology. 
Now the assertion follows from Proposition \ref{strict}.
\end{proof}

\begin{proposition}\label{tag3}
The element $U_T\in\cM(C^*(T))$ defined in Proposition \ref{tga2}
satisfies the following.
\benu
\item $U_T^*U_T=1$.
\item $T^0(f)U_T=U_TT^0(\hat{\sigma}(f))$ for $f\in C_0(X)$. 
\item $T^n(f)=U_T^nT^0(f)$ for $f\in C_0(X)$
\eenu
\end{proposition}

\begin{proof}
\benu
\item For an approximate unit $u_i$ of $C_0(X)$, 
$|u_i|^2$ is also an approximate unit. 
Hence $T^1(u_i)^*T^1(u_i)=T^0(|u_i|^2)$ 
converges strictly to $1\in\cM(C^*(T))$. 
Thus we have $U_T^*U_T=1$.
\item By Proposition \ref{tga2}, 
$T^0(f)U_T=T^1(\hat{\sigma}(f))=U_TT^0(\hat{\sigma}(f))$. 
\item Similarly as the proof of Proposition \ref{tga2},
we have $U_TT^{n-1}(f)=T^n(f)$. 
Hence we get $T^n(f)=U_T^nT^0(f)$. 
\eenu
\end{proof}

Thus from a Toeplitz $E_\Sigma$-pair $T$, 
we get an isometry $U_T\in\cM(C^*(T))$ 
satisfying $T^0(f)U_T=U_TT^0(\hat{\sigma}(f))$. 
Conversely, we have the following. 

\begin{proposition}\label{tag4}
Let $A$ be a $C^*$-algebra. 
Suppose that a $*$-homomorphism $\pi: C_0(X)\to A$ 
and an isometry $u\in\cM(A)$ 
satisfy that $\pi(f)u=u\pi(\hat{\sigma}(f))$ for all $f\in C_0(X)$. 
Define $T^0,T^1\colon C_0(X)\to A$ by $T^0=\pi$ and $T^1(f)=u\pi(f)$. 
Then $T=(T^0,T^1)$ is a Toeplitz $E_\Sigma$-pair. 
\end{proposition}

\begin{proof}
Since $u$ is an isometry, 
$T$ satisfies the condition (i) in Definition \ref{DefTpl}.
The condition (ii) is easily checked 
from the relation $\pi(f)u=u\pi(\hat{\sigma}(f))$. 
Thus $T=(T^0,T^1)$ is a Toeplitz $E_\Sigma$-pair. 
\end{proof}

A pair $\{\pi,u\}$ appeared in Proposition \ref{tag4} 
can be considered as a kind of covariant representations of $(X,\sigma)$. 
These representations correspond to Toeplitz $E_\Sigma$-pairs. 
We study which representation $\{\pi,u\}$ 
corresponds to a Cuntz-Krieger $E_\Sigma$-pair. 
Since $\sigma$ is proper, we have $(E_\Sigma^0)_{\rs{fin}}=X$ 
and $(E_\Sigma^0)_{\rs{sce}}=X\setminus\sigma(X)$. 
Hence $(E_\Sigma^0)_{\rs{rg}}=X\setminus\overline{X\setminus\sigma(X)}$ 
which is the interior of the image $\sigma(X)$ of $\sigma$. 
It is not difficult to see that 
the map $\varPhi^1\colon \cK(C_d(E_\Sigma^1))\to C^*(T)$ 
can be expressed as $\varPhi^1(f)=U_TT^0(f)U_T^*$ for $f\in C_0(X)$ 
by identifying $\cK(C_d(E_\Sigma^1))$ with $C_0(X)$. 
Hence we have $\varPhi^1(\pi_r(f))=U_TT^0(\hat{\sigma}(f))U_T^*$ 
for $f\in C_0(X)$. 
We write $P_T=1-U_TU_T^*\in\cM(C^*(T))$ which is a projection. 
For $f\in C_0(X)$, 
we have 
$$T^0(f)P_T=T^0(f)-U_TT^0(\hat{\sigma}(f))U_T^*=P_TT^0(f)$$
by Proposition \ref{tag3} (ii). 
By the above argument, we have the following. 

\begin{proposition}\label{tag5}
For a Toeplitz $E_\Sigma$-pair $T$, 
the following are equivalent.
\benu
\item $T$ is a Cuntz-Krieger $E_\Sigma$-pair.
\item $T^0(f)=U_TT^0(\hat{\sigma}(f))U_T^*$
for all $f\in C_0(X)$ with $f(x)=0$ for $x\notin\sigma(X)$.
\item $T^0(f)P_T=0$ 
for all $f\in C_0(X)$ with $f(x)=0$ for $x\notin\sigma(X)$.
\eenu
\end{proposition}

\begin{proposition}\label{tag6}
The $C^*$-algebra $\cO(E_\Sigma)$ is the universal $C^*$-algebra 
generated by products of a copy of $C_0(X)$ 
and an isometry $u$ satisfying that 
\benu
\item $fu=u\hat{\sigma}(f)$ for $f\in C_0(X)$,  
\item $f=u\hat{\sigma}(f)u^*$ for $f\in C_0(X)$ 
with $f(x)=0$ for $x\notin\sigma(X)$. 
\eenu
\end{proposition}

When $\sigma$ is surjective, 
for a Cuntz-Krieger $E_\Sigma$-pair $T$ 
we have $P_T=0$, that is, $U_T$ is a unitary. 
Hence, we have the following. 

\begin{corollary}
When $\sigma$ is surjective, 
the $C^*$-algebra $\cO(E_\Sigma)$ is the universal $C^*$-algebra 
generated by products of a copy of $C_0(X)$ 
and a unitary $u$ satisfying that 
$\hat{\sigma}(f)=u^*fu$ for $f\in C_0(X)$. 
\end{corollary}

\begin{corollary}\label{tag7}
When $\sigma$ is a homeomorphism, 
Cuntz-Krieger $E_\Sigma$-pairs correspond to 
covariant representations of the dynamical system $\Sigma$, 
and we have a natural isomorphism 
between $\cO(E_\Sigma)$ and 
the homeomorphism $C^*$-algebra $A(\Sigma)$. 
\end{corollary}

It is complicated to describe 
Toeplitz $E_\Sigma$-pairs or Cuntz-Krieger $E_\Sigma$-pairs 
when $\sigma$ is not proper, or 
when $\sigma$ is defined just on some open subset of $X$.

\section{Fock representations}\label{Fock}

The purpose of this section is 
a construction of an injective Cuntz-Krieger $E$-pair $\tau=(\tau^0,\tau^1)$
of a graph $E$ by using the so-called Fock space $C_d(E^{*})$.
The map $\rho_\tau\colon \cO(E)\to C^*(\tau)$
is called the Fock representation. 
In the next section, 
it will turn out that the Fock representation is faithful. 
Hence the construction done in this section gives us 
a concrete description of $\cO(E)$. 
The results here will be used in Section \ref{Kgroups}
to compute $K$-groups. 

\begin{definition}
For a topological graph $E=(E^0,E^1,d,r)$, 
we denote by $E^{*}$ the disjoint union of 
$E^0,E^1,\ldots,E^n,\ldots$. 
\end{definition}

The set $E^{*}$ is called a {\em finite path space} 
of a topological graph $E$. 
We can define $d,r\colon E^{*}\to E^0$ by using $d^n,r^n$. 
The $C^*$-correspondence $C_d(E^{*})$ over $C_0(E^0)$, 
which is called a {\em Fock space}, 
is isomorphic to $\bigoplus_{n=0}^\infty C_{d^n}(E^n)$. 
We denote the left action of $C_0(E^0)$ on $C_d(E^{*})$ 
by $\sigma^0\colon C_0(E^0)\to\cL(C_d(E^{*}))$. 
Explicitly, for $f\in C_0(E^0)$ we see that 
$$\begin{cases}
\sigma^0(f)\xi=f \xi, &\mbox{for }\xi\in C_{d^0}(E^0)\subset C_d(E^{*}),\\
\sigma^0(f)(\xi\otimes\eta)=(\pi_r(f)\xi)\otimes\eta,
&\mbox{for }\xi\in C_{d^1}(E^1),\eta\in C_{d^n}(E^n)\ (n\in\N),\\
\end{cases}$$
where $\pi_r$ is the left action of $C_0(E^0)$ on $C_d(E^1)$
defined in Section \ref{graphalg}.
We define a linear map 
$\sigma^1\colon C_d(E^1)\ni\xi\mapsto\sigma^1(\xi)\in\cL(C_d(E^{*}))$ by
$\sigma^1(\xi)\eta=\xi\otimes\eta\in C_{d^{n+1}}(E^{n+1})$ 
for $\eta\in C_{d^n}(E^n)\subset C_d(E^{*})$.
The routine computation shows the following formulae 
of the adjoint $\sigma^1(\xi_0)^*$ of $\sigma^1(\xi_0)$ 
for $\xi_0\in C_d(E^1)$: 
$$\begin{cases}
\sigma^1(\xi_0)^*\xi=0, &\mbox{for }\xi\in C_{d^0}(E^0),\\
\sigma^1(\xi_0)^*(\xi\otimes\eta)=\sigma^0(\ip{\xi_0}{\xi})\eta,
&\mbox{for }\xi\in C_{d^1}(E^1),\eta\in C_{d^n}(E^n)\ (n\in\N).\\
\end{cases}$$
Now, it is easy to see the following.

\begin{proposition}[cf. {\cite[Proposition 1.3]{Pi}}]
The pair $\sigma=(\sigma^0,\sigma^1)$ is a Toeplitz $E$-pair.
\end{proposition}

Recall that the map 
$\varPhi^1\colon \cK(C_d(E^1))\to C^*(\sigma)\subset\cL(C_d(E^{*}))$ 
is defined by $\varPhi^1(\theta_{\xi,\eta})=\sigma^1(\xi)\sigma^1(\eta)^*$.
For $x\in \cK(C_d(E^1))$,
we see that 
$$\begin{cases}
\varPhi^1(x)\xi=0, &\mbox{for }\xi\in C_{d^0}(E^0),\\
\varPhi^1(x)(\xi\otimes\eta)=
(x\xi)\otimes\eta,
&\mbox{for }\xi\in C_{d^1}(E^1),\eta\in C_{d^n}(E^n)\ (n\in\N).\\
\end{cases}$$
From this computation, 
we get the following lemma, which measures 
how far the Toeplitz $E$-pair $\sigma=(\sigma^0,\sigma^1)$ is 
from being a Cuntz-Krieger $E$-pair. 

\begin{lemma}\label{differnce}
For $f\in C_0(\s{E}{0}{rg})$,
take $\xi_0,\eta_0\in C_{d^0}(E^0)\subset C_d(E^{*})$ 
with $\xi_0\overline{\eta_0}=f$.
Then we have $\theta_{\xi_0,\eta_0}=\sigma^0(f)-\varPhi^1(\pi_r(f))$.
\end{lemma}

\begin{proof}
By the computation above, 
it suffices to show that 
$\theta_{\xi_0,\eta_0}(\xi)=f\xi$
if $\xi\in C_{d^0}(E^0)\subset C_d(E^{*})$, and 
$\theta_{\xi_0,\eta_0}(\xi)=0$ 
if $\xi\in C_{d^n}(E^n)\subset C_d(E^{*})$ for $n\geq 1$. 
The former is verified by 
$$\theta_{\xi_0,\eta_0}(\xi)=\xi_0\overline{\eta_0}\xi=f\xi,$$
and the latter is obvious. 
The proof is completed. 
\end{proof}

For each $n\in\N$,
we define an open subset $\s{E}{n}{rg}$ of $E^n$ 
by $\s{E}{n}{rg}=(d^n)^{-1}(\s{E}{0}{rg})$ 
and an open subset $\s{E}{*}{rg}$ of $E^{*}$ 
by $\s{E}{*}{rg}=d^{-1}(\s{E}{0}{rg})$. 
Note that 
$C_d(\s{E}{*}{rg})
=\bigoplus_{n=0}^\infty C_{d^n}(\s{E}{n}{rg})\subset C_d(E^{*})$ 
and that $\cK(C_d(\s{E}{*}{rg}))$ is an ideal of $\cL(C_d(E^{*}))$.

\begin{proposition}
We have $\cK(C_d(\s{E}{*}{rg}))\subset C^*(\sigma)$.
\end{proposition}

\begin{proof}
It suffices to show that $\theta_{\xi,\eta}\in C^*(\sigma)$ 
for $\xi\in C_{d^n}(\s{E}{n}{rg})\subset C_d(\s{E}{*}{rg})$ and 
$\eta\in C_{d^m}(\s{E}{m}{rg})\subset C_d(\s{E}{*}{rg})$ 
for $n,m\in\N$.
By Lemma \ref{surjection2}, 
we can find 
$\xi'\in C_{d^n}(E^n)$, $\eta'\in C_{d^m}(E^m)$ and $f,g\in C_0(\s{E}{0}{rg})$
with $\xi=\xi'f$ and $\eta=\eta'g$.
Once we consider $f,g$ as elements of $C_{d^0}(E^0)\subset C_d(E^{*})$,
we have $\xi=\sigma^n(\xi')f$ and $\eta=\sigma^m(\eta')g$. 
By Lemma \ref{differnce}, 
we have that 
$$\theta_{f,g}=\sigma^0(h)-\varPhi^1(\pi_r(h))\in C^*(\sigma),$$
where $h=f\overline{g}\in C_0(\s{E}{0}{rg})$.
Hence we get
$$\theta_{\xi,\eta}=\sigma^n(\xi')\theta_{f,g}\sigma^m(\eta')^*
\in C^*(\sigma).$$
The proof is completed.
\end{proof}

Let $\tau^0\colon C_0(E^0)\to\cL(C_d(E^{*}))/\cK(C_d(\s{E}{*}{rg}))$ 
and $\tau^1\colon C_d(E^1)\to\cL(C_d(E^{*}))/\cK(C_d(\s{E}{*}{rg}))$ 
be the compositions 
of the natural surjection 
$\cL(C_d(E^{*}))\to \cL(C_d(E^{*}))/\cK(C_d(\s{E}{*}{rg}))$
with $\sigma^0$ and $\sigma^1$ respectively.
By Lemma \ref{differnce}, 
the pair $\tau=(\tau^0,\tau^1)$ is a Cuntz-Krieger $E$-pair.
We will show that this pair is injective.

\begin{lemma}\label{leminj}
Let $n$ be a natural number,
and $f$ be an element of $C_0(E^0)$ 
with $\pi_{r^n}(f)\in \cK(C_{d^n}(\s{E}{n}{rg}))$.
If $e\in E^n$ satisfies $f(r^n(e))\neq 0$, 
then $e\in  \s{E}{n}{rg}$.
\end{lemma}

\begin{proof}
To the contrary, assume that $e\in E^n$ satisfies 
$f(r^n(e))\neq 0$ and $e\notin  \s{E}{n}{rg}$.
Take $\xi_1,\ldots,\xi_m$, $\eta_1,\ldots,\eta_m$ in $C_c(\s{E}{n}{rg})$
arbitrarily,
and we will show that 
$$\bigg\|\pi_r(f)-\sum_{k=1}^m\theta_{\xi_k,\eta_k}\bigg\|\geq |f(r^n(e))|,$$
which contradicts the fact that $\pi_{r^n}(f)\in \cK(C_{d^n}(\s{E}{n}{rg}))$.
Set $X=\bigcup_{k=1}^m\supp(\eta_k)$ which is a compact set 
with $X\subset  \s{E}{n}{rg}$.
We can find a neighborhood $U$ of $e$ such that $X\cap U=\emptyset$
and that the restriction of $d^n$ to $U$ is injective.
Take $\zeta\in C_c(U)\subset C_{d^n}(E^n)$ 
with $0\leq\zeta\leq 1$ and $\zeta(e)=1$.
Similarly as the proof of Proposition \ref{piK},
we have that $\|\zeta\|=1$ and 
$$\bigg\|\big(\pi_{r^n}(f)-\sum_{k=1}^m\theta_{\xi_k,\eta_k}\big)\zeta\bigg\|=
\|\pi_{r^n}(f)\zeta\|\geq |f(r^n(e))|.$$
We are done.
\end{proof}

\begin{proposition}\label{FockRep}
The Cuntz-Krieger $E$-pair $\tau=(\tau^0,\tau^1)$ is injective.
\end{proposition}

\begin{proof}
To the contrary, assume that there exists $f\in C_0(E^0)$ 
with $f\neq 0$ and $\sigma^0(f)\in \cK(C_d(\s{E}{*}{rg}))$.
The fact $\sigma^0(f)\in \cK(C_d(\s{E}{*}{rg}))$ implies that 
$\pi_{r^n}(f)\in \cK(C_{d^n}(\s{E}{n}{rg}))$ for every $n\in\N$. 
We can find $\e>0$ and a non-empty open subset $V$ of $E^0$ such that 
$|f(v)|\geq\e$ for any $v\in V$. 
We will show that $(r^n)^{-1}(V)$ is a non-empty subset of $\s{E}{n}{rg}$ 
for every $n\in\N$ by induction. 
For $n=0$, we have $V\subset \s{E}{0}{rg}$ 
because $\pi_{r^0}(f)\in \cK(C_{d^0}(\s{E}{0}{rg}))\cong C_0(\s{E}{0}{rg})$. 
Assume that $(r^n)^{-1}(V)$ is a non-empty subset of $\s{E}{n}{rg}$. 
Then we have 
$d^n((r^n)^{-1}(V))\subset \s{E}{0}{rg}\subset \overline{r(E^1)}$. 
Since $d^n((r^n)^{-1}(V))$ is non-empty and open, 
there exists $e\in (r^n)^{-1}(V)$ such that $d^n(e)\in r(E^1)$. 
Hence $(r^{n+1})^{-1}(V)$ is non-empty. 
Since $|f(r^{n+1}(e))|\geq\e$ for $e\in (r^{n+1})^{-1}(V)$, 
we have $(r^{n+1})^{-1}(V)\subset\s{E}{n+1}{rg}$ by Lemma \ref{leminj}. 
Thus we have shown that $(r^n)^{-1}(V)$ is a non-empty subset 
of $\s{E}{n}{rg}$ for every $n\in\N$. 
We will show that $\|\sigma^0(f)-\sum_{k=1}^m\theta_{\xi_k,\eta_k}\|\geq\e$ 
for any $\xi_k,\eta_k\in C_d(\s{E}{*}{rg})$, 
which contradicts the fact that $\sigma^0(f)\in \cK(C_d(\s{E}{*}{rg}))$. 
To this end, it suffices to find $\zeta_n\in C_{d^n}(E^n)$ 
with $\|\zeta_n\|=1$, $\|\sigma^0(f)\zeta_n\|\geq\e$ for each $n\in\N$. 
Since $(r^n)^{-1}(V)$ is not empty, 
we can find $\zeta_n\in C_c((r^n)^{-1}(V))\subset C_{d^n}(E^n)$ 
with $\|\zeta_n\|=1$. 
Since 
$$|(\sigma^0(f)\zeta_n)(e)|=|f(r^n(e))\ \zeta_n(e)|\geq\e|\zeta_n(e)|$$
for $e\in (r^n)^{-1}(V)$, 
we have $\|\sigma^0(f)\zeta_n\|\geq\e\|\zeta_n\|=\e$. 
We are done. 
\end{proof}

As claimed in the previous section, 
Proposition \ref{FockRep} implies the following. 

\begin{proposition}\label{univinj}
The universal Cuntz-Krieger $E$-pair $t=(t^0,t^1)$ is injective. 
\end{proposition}

The map $\rho_{\tau}\colon \cO(E)\to C^*(\tau)$ 
is called the {\em Fock representation}.
In the next section, 
we will show that the Fock representation gives us an isomorphism
$C^*(\tau)\cong\cO(E)$ (Corollary \ref{Fockrep}).

We finish this section by stating a relation 
between our $C^*$-algebras $\cT(E)$, $\cO(E)$ 
and ones defined in \cite{Pi}. 
The $C^*$-algebra $C^*(\sigma)\subset\cL(C_d(E^{*}))$ is 
exactly the same as 
the augmented Toeplitz algebra $\widetilde{\cT}_{C_d(E^1)}$ 
of the $C^*$-correspondence $C_d(E^1)$ over $C_0(E^0)$ 
defined in \cite[Remark 1.2 (3)]{Pi}.
Hence Theorem 3.4 of \cite{Pi} gives the following 
because the conditions there are 
the same as the ones of Toeplitz $E$-pairs.

\begin{proposition}\label{univToep}
The pair $\sigma=(\sigma^0,\sigma^1)$ is the universal 
Toeplitz $E$-pair. 
Hence $C^*(\sigma)$ is isomorphic to $\cT(E)$.
\end{proposition}

Therefore we see that $\cT(E)$ is isomorphic to 
the augmented Toeplitz algebra $\widetilde\cT_{C_d(E^1)}$. 
Of course one can show Proposition \ref{univToep} 
by using a similar argument in Section \ref{secGIUT} 
(which is actually the same as the proof in \cite{Pi}). 
One can also deduce Proposition \ref{univToep} from Theorem \ref{GIUT}, 
by proving that Toeplitz pairs corresponds bijectively to 
Cuntz-Krieger pairs of a certain topological graph 
(see \cite{Ka2}). 

The augmented Cuntz-Pimsner algebra $\widetilde\cO_{C_d(E^1)}$ 
is isomorphic to the universal $C^*$-algebra generated 
by a Toeplitz $E$-pair $T=(T^0,T^1)$ satisfying
$T^0(f)=\varPhi^1(\pi_r(f))$ 
for every $f\in C_0(\s{E}{0}{fin})$ not only $f\in C_0(\s{E}{0}{rg})$
(see, \cite[Theorem 3.12]{Pi}).
Hence there exists a surjection $\cO(E)\to\widetilde\cO_{C_d(E^1)}$.
In this sense, 
$\widetilde\cO_{C_d(E^1)}$ is ``smaller'' than $\cO(E)$.
Sometimes $\widetilde\cO_{C_d(E^1)}$ is too small and can be $0$.
If there exists a source, 
then $t^0\colon C_0(E^0)\to\widetilde\cO_{C_d(E^1)}$ never become injective
because $\pi_r(f)=0$ for $f\in C_0(\s{E}{0}{sce})$ 
by Proposition \ref{pir-1}.
In the case that there exist no sources, we have the following.

\begin{proposition}\label{CP}
If $\overline{r(E^1)}=E^0$, 
then $\cO(E)\cong\widetilde\cO_{C_d(E^1)}$. 
\end{proposition}

\begin{remark}
In \cite{D}, 
V. Deaconu introduced compact graphs 
which are topological graphs $E=(E^0,E^1,d,r)$ 
such that both $E^0$ and $E^1$ are compact, 
and both $d$ and $r$ are 
surjective and locally homeomorphic. 
He associated a $C^*$-algebra with a compact graph 
by constructing a certain groupoid, 
and showed in \cite[Theorem 4.3]{D} 
that this is isomorphic to 
the Cuntz-Pimsner algebras of 
the $C^*$-correspondence defined by the compact graph. 
Hence Proposition \ref{CP} 
implies that his $C^*$-algebras are isomorphic 
to our $C^*$-algebras. 
As he pointed out in the last part of \cite{D}, 
the $C^*$-algebras arising from polymorphisms defined in \cite{AR} 
are different from our $C^*$-algebras in general.
\end{remark}

\begin{remark}
We should note that the $C^*$-algebra $C^*(\tau)$ is the same as 
the relative Cuntz-Pimsner algebra $\cO(C_0(\s{E}{0}{rg}),C_d(E^1))$ 
determined by the ideal $C_0(\s{E}{0}{rg})\subset C_0(E^0)$ 
\cite[Definition 2.18]{MS}, 
and Proposition \ref{FockRep} follows from \cite[Proposition 2.21]{MS}. 
By Corollary \ref{Fockrep}, 
we have the isomorphism 
$\cO(E)\cong \cO(C_0(\s{E}{0}{rg}),C_d(E^1))$ 
(this also can be proved using \cite[Theorem 2.19]{MS}).
\end{remark}

\section{The gauge-invariant uniqueness theorem}\label{secGIUT}

In this section and the next section,
we investigate for which Cuntz-Krieger $E$-pair $T$, 
$\rho_T$ gives an isomorphism $C^*(T)\cong \cO(E)$.
A pair is necessarily injective, 
but this condition is not sufficient in general.
In this section, 
we give two kinds of extra conditions 
for the isomorphism $C^*(T)\cong \cO(E)$,
namely the existence of gauge actions 
and the existence of conditional expectations (Theorem \ref{GIUT}).
In the next section,
we deal with the problem of determining topological graphs for which 
injectivity of $T$ is sufficient for the isomorphism $C^*(T)\cong \cO(E)$.

By the universality of $\cO(E)$, 
there exists an action $\beta\colon \T\curvearrowright \cO(E)$ defined by 
$\beta_z(t^0(f))=t^0(f)$ and $\beta_z(t^1(\xi))=zt^1(\xi)$ 
for $f\in C_0(E^0)$, $\xi\in C_d(E^1)$ and $z\in\T$.
The action $\beta$ is called the {\em gauge action}. 
It is easy to see that 
$$\beta_z(t^n(\xi)t^m(\eta)^*)=z^{n-m}t^n(\xi)t^m(\eta)^*\qquad 
( \xi\in C_d(E^n),\ \eta\in C_d(E^m) ).$$
We define a linear map $\varPsi\colon \cO(E)\to \cO(E)$ by
$$\varPsi(x)=\int_{\T}\beta_z(x)dz$$
where $dz$ is the normalized Haar measure of $\T$.
Then we have 
$$\varPsi(t^n(\xi)t^m(\eta)^*)=\delta_{n,m}t^n(\xi)t^m(\eta)^*\qquad 
( \xi\in C_d(E^n),\ \eta\in C_d(E^m) ),$$
where $\delta_{n,m}$ is the Kronecker delta.
Hence $\varPsi$ is a faithful conditional expectation onto a subalgebra
$$\F:=\cspa\{t^k(\xi)t^k(\eta)^*\mid \xi,\eta\in C_d(E^k),\ k\in\N\}.$$

\begin{definition}
For a Cuntz-Krieger $E$-pair $T=(T^0,T^1)$, 
we define subalgebras $\F_T^n,\G_T^n$ for $n\in\N$ and $\F_T$ of $C^*(T)$ by
\begin{align*}
\G_T^{n}
&=\cspa\{T^n(\xi)T^n(\eta)^*\mid \xi,\eta\in C_d(E^n)\},\\
\F_T^{n}
&=\cspa\{T^k(\xi)T^k(\eta)^*\mid \xi,\eta\in C_d(E^k),\ 0\leq k\leq n\},\\
\F_T&=\cspa\{T^k(\xi)T^k(\eta)^*\mid \xi,\eta\in C_d(E^k),\ k\in\N\}.
\end{align*}
We simply write $\F^{n},\G^{n}$ for $n\in\N$ and $\F$ 
for the corresponding subalgebras in $\cO(E)$.
\end{definition}

Note that $\G_T^{n}$ is an ideal of the $C^*$-algebra $\F_T^{n}$ 
and that $\F_T^{n+1}=\G_T^{n+1}+\F_T^{n}$ for each $n\in\N$. 
We also see that $\F_T=\overline{\bigcup_{n=0}^\infty \F_T^{n}}$. 
Note that $\G_T^{n}$ is the image of $\varPhi^n$,
hence $\cK(C_d(E^n))\cong\G_T^{n}$ when $T$ is injective.
We will show that if a Cuntz-Krieger $E$-pair $T$ is injective, 
then the restriction of $\rho_T$ to $\F$ is 
an isomorphism onto $\F_T$ (Proposition \ref{isomF}).

\begin{lemma}\label{FG1}
For an injective Cuntz-Krieger $E$-pair $T$, 
we have $\G_T^{0}\cap \G_T^{1}=T^0(C_0(\s{E}{0}{rg}))$.
\end{lemma}

\begin{proof}
This follows from the definition of Cuntz-Krieger pairs and 
Proposition \ref{T=Pp}. 
\end{proof}

\begin{lemma}\label{FG2}
If a Cuntz-Krieger $E$-pair $T$ is injective, 
then $\F_T^{n}\cap \G_T^{n+1}=\G_T^{n}\cap \G_T^{n+1}
=\varPhi^n\big(\cK(C_d(\s{E}{n}{rg}))\big)$
for every $n\in\N$.
\end{lemma}

\begin{proof}
Take $x\in\F_T^{n}\cap \G_T^{n+1}$. 
Let $\{u_i\}_{i\in\I}$ be an approximate unit of $\G_T^{n}$. 
Since 
$$\G_T^{n+1}=\cspa\{T^n(\xi)T^1(\xi')T^1(\eta')^*T^n(\eta)^*\mid 
\xi,\eta\in C_d(E^n),\ \xi',\eta'\in C_d(E^1)\},$$
$\{u_i\}_{i\in\I}$ is also an approximate unit of $\G_T^{n+1}$.
Hence we have $x=\lim_i u_ix$.
Since $x\in\F_T^{n}$ and $\G_T^{n}$ is an ideal of $\F_T^{n}$, 
we have $u_ix\in\G_T^{n}$ for $i\in\I$.
Hence $x\in\G_T^{n}$. 
Thus we have $\F_T^{n}\cap \G_T^{n+1}=\G_T^{n}\cap \G_T^{n+1}$.

Take $\xi,\eta\in C_d(\s{E}{n}{rg})$ arbitrarily.
We can find $\xi'\in C_d(E^n)$ and $f\in C_0(\s{E}{0}{rg})$ with $\xi=\xi' f$.
Since $T$ is a Cuntz-Krieger $E$-pair, 
we have 
$$\varPhi^n(\theta_{\xi,\eta})=
T^n(\xi)T^n(\eta)^*=T^n(\xi')T^0(f)T^n(\eta)^*
=T^n(\xi')\varPhi^1(\pi_r(f))T^n(\eta)^*\in\G_T^{n+1}.$$
Thus $\varPhi^n\big(\cK(C_d(\s{E}{n}{rg}))\big)\subset \G_T^{n}\cap \G_T^{n+1}$.
Conversely take $x\in\cK(C_d(E^n))$ 
with $\varPhi^n(x)\in \G_T^{n}\cap \G_T^{n+1}$.
For $\xi,\eta\in C_d(E^n)$, 
we have 
$$T^0(\ip{\xi}{x\eta})
=T^n(\xi)^*\varPhi^n(x)T^n(\eta)\in\G_T^{0}\cap \G_T^{1}.$$
By Lemma \ref{FG1}, we have $\ip{\xi}{x\eta}\in C_0(\s{E}{0}{rg})$.
Therefore we get $x\in \cK(C_d(\s{E}{n}{rg}))$ 
by Lemma \ref{surjection4}.
Thus we have shown that 
$\F_T^{n}\cap \G_T^{n+1}=\G_T^{n}\cap \G_T^{n+1}
=\varPhi^n\big(\cK(C_d(\s{E}{n}{rg}))\big)$.
\end{proof}

\begin{proposition}\label{isomF}
For an injective Cuntz-Krieger $E$-pair $T$, 
the restriction of $\rho_T$ to $\F$ is 
an isomorphism onto $\F_T$.
\end{proposition}

\begin{proof}
First note that the restriction of $\rho_T$ to $\G^{n}$ is 
an isomorphism onto $\G_T^{n}$ for every $n\in\N$ 
because $\rho_T\circ\varphi^n=\varPhi^n$ 
and $\varphi^n,\varPhi^n$ are isomorphisms 
onto $\G^{n}$ and $\G_T^{n}$ respectively.
To finish the proof, it suffices to show that 
the restriction of $\rho_T$ to $\F^{n}$ is 
an isomorphism onto $\F_T^{n}$ for every $n\in\N$.
We will prove this by induction with respect to $n\in\N$.
The restriction of $\rho_T$ to $\F^{0}$ is 
an isomorphism onto $\F_T^{0}$ 
because $\F^{0}=\G^{0}$ and $\F_T^{0}=\G_T^{0}$.
Assume that the restriction of $\rho_T$ to $\F^{n}$ is 
an isomorphism onto $\F_T^{n}$.
We have the following commutative diagram with exact row
$$\begin{CD}
0 @>>> \G^{n+1} @>>> \F^{n+1} @>>> \F^{n+1}/\G^{n+1} @>>> 0\phantom{.} \\
@. @VV\rho_TV @VV\rho_TV  @VVV \\
0 @>>> \G_T^{n+1} @>>> \F_T^{n+1} @>>> \F_T^{n+1}/\G_T^{n+1} @>>> 0. \\
\end{CD}$$
To prove that the restriction of $\rho_T$ to $\F^{n+1}$ is 
an isomorphism onto $\F_T^{n+1}$,
it is sufficient to see that the map 
$\F^{n+1}/\G^{n+1}\to\F_T^{n+1}/\G_T^{n+1}$ induced by $\rho_T$ 
is an isomorphism. 
Since $\F_T^{n+1}=\F_T^{n}+\G_T^{n+1}$,
we have 
$$\F_T^{n+1}/\G_T^{n+1}\cong \F_T^{n}/(\F_T^{n}\cap\G_T^{n+1})
=\F_T^{n}/(\G_T^{n}\cap\G_T^{n+1}).$$
By the assumption of the induction, 
the restriction of $\rho_T$ to $\F^{n}$ is 
an isomorphism onto $\F_T^{n}$. 
By Lemma \ref{FG2}, we have
$$\G^{n}\cap\G^{n+1}\cong\cK(C_d(\s{E}{n}{rg}))\cong \G_T^{n}\cap\G_T^{n+1}.$$
Hence the restriction of $\rho_T$ to $\G^{n}\cap\G^{n+1}$ 
is an isomorphism onto $\G_T^{n}\cap\G_T^{n+1}$. 
Therefore the map 
$\F^{n}/(\G^{n}\cap\G^{n+1})\to\F_T^{n}/(\G_T^{n}\cap\G_T^{n+1})$ 
induced by $\rho_T$ is an isomorphism. 
Hence the map 
$\F^{n+1}/\G^{n+1}\to\F_T^{n+1}/\G_T^{n+1}$ induced by $\rho_T$ 
is also an isomorphism. 
Thus we have shown that the restriction of $\rho_T$ to $\F^{n+1}$ is 
an isomorphism onto $\F_T^{n+1}$. 
We are done.
\end{proof}

Now we have the following gauge-invariant uniqueness theorem.

\begin{theorem}\label{GIUT}
For a topological graph $E=(E^0,E^1,d,r)$ and 
a Cuntz-Krieger $E$-pair $T=(T^0,T^1)$, 
the following are equivalent:
\benu
\item The map $\rho_{T}\colon \cO(E)\to C^*(T)$ is an isomorphism.
\item The map $T^0$ is injective and 
there exists an automorphism $\beta'_z$ of $C^*(T)$ 
such that $\beta'_z(T^0(f))=T^0(f)$ and $\beta'_z(T^1(\xi))=zT^1(\xi)$
for every $z\in\T$.
\item The map $T^0$ is injective and 
there exists a conditional expectation $\varPsi_T$ 
from $C^*(T)$ onto $\F_T$ such that 
$\varPsi_T(T^n(\xi)T^m(\eta)^*)=\delta_{n,m}T^n(\xi)T^m(\eta)^*$ 
for $\xi\in C_d(E^n)$ and $\eta\in C_d(E^m)$.
\eenu
\end{theorem}

\begin{proof}
(i)$\Rightarrow$(ii): Already shown.

(ii)$\Rightarrow$(iii): Set $\varPsi_T(x)=\int_{\T}\beta'_z(x)dz$.

(iii)$\Rightarrow$(i): Since the map $T^0$ is injective, 
we see that the restriction of $\rho_T$ to $\F$ is 
an isomorphism onto $\F_T$ by Proposition \ref{isomF}. 
Now we see that the map $\rho_{T}\colon \cO(E)\to C^*(T)$ is an isomorphism
by the standard argument of conditional expectations 
(see, for example, \cite[Proposition 3.11]{Ka1}). 
\end{proof}

\begin{remark}
If there exists an automorphism $\beta'_z$ of $C^*(T)$ 
such that $\beta'_z(T^0(f))=T^0(f)$ and $\beta'_z(T^1(\xi))=zT^1(\xi)$
for every $z\in\T$, 
then $\beta'\colon \T\ni z\mapsto \beta'_z\in\Aut(C^*(T))$ 
becomes automatically a strongly continuous homomorphism. 
This fact is used implicitly in the proof of the implication 
(ii)$\Rightarrow$(iii) in the above theorem. 
\end{remark}

\begin{corollary}\label{Fockrep}
The Fock representation $\rho_\tau$ of $\cO(E)$ is faithful.
\end{corollary}

\begin{proof}
We check the condition (ii) in Theorem \ref{GIUT}.
We have already seen that $\tau^0$ is injective in Proposition \ref{FockRep}.
For $z\in\T$, 
define a unitary $u_z\in\cL(C_d(E^{*}))$ by $u_z(\xi)=z^n\xi$ 
for $\xi\in C_{d^n}(E^n)\subset C_d(E^{*})$ ($n\in\N$). 
We define an automorphism $\beta_z'=\Ad(\pi(u_z))$
of $\cL(C_d(E^{*}))/\cK(C_d(\s{E}{*}{rg}))$ 
by $\beta_z'(x)=\pi(u_z)x\pi(u_z)^*$, 
where $\pi$ is the natural surjection 
$$\cL(C_d(E^{*}))\to\cL(C_d(E^{*}))/\cK(C_d(\s{E}{*}{rg})).$$
One can easily see that 
$\beta'_z(\tau^0(f))=\tau^0(f)$ and $\beta'_z(\tau^1(\xi))=z\tau^1(\xi)$
for $f\in C_0(E^0)$ and $\xi\in C_d(E^1)$.
This implies that $\beta_z'$ fixes $C^*(\tau)$ globally for each $z\in\T$.
Hence the restriction of $\beta_z'$ to $C^*(\tau)$ 
is an automorphism of $C^*(\tau)$. 
By Theorem \ref{GIUT}, 
we see that $\rho_\tau\colon \cO(E)\to C^*(\tau)$ 
is an isomorphism.
\end{proof}

\section{The Cuntz-Krieger uniqueness theorem}\label{secCKUT}

In this section, we see that 
if a graph $E$ satisfies a certain condition, 
then the condition that $T^0$ is injective is 
not only necessary but sufficient for $C^*(T)\cong \cO(E)$. 
To this end, we need a more precise description of $\F_T^{n}$
for an injective Cuntz-Krieger $E$-pair $T$. 
Let us fix an injective Cuntz-Krieger $E$-pair $T$ for a while. 

For each $n\in\N$,
we define a $*$-homomorphism $\pi_n^{n}\colon \F_T^{n}\to\cL(C_d(E^n))$ by 
$$T^n(\pi_n^{n}(x)\xi)=xT^n(\xi)\quad\mbox{ for }
x\in\F_T^{n},\ \xi\in C_d(E^n).$$
Note that $T^n\colon C_d(E^n)\to C^*(T)$ is injective 
and that $xT^n(\xi)$ lies in the image of $T^n$ 
for $x\in\F_T^{n},\ \xi\in C_d(E^n)$. 
The restriction of $\pi_n^n$ to $\G_T^n$ coincides 
with the isomorphism $\varPhi^n\colon \G_T^n\to\cK(C_d(E^n))$. 

We define a closed subset $\s{E}{n}{sg}$ of $E^n$ by 
$\s{E}{n}{sg}=(d^n)^{-1}(\s{E}{0}{sg})=E^n\setminus \s{E}{n}{rg}$ 
for each $n\in\N$.
Recall that there exists a $*$-homomorphism 
$\omega^n\colon \cL(C_d(E^n))\to\cL(C_d(\s{E}{n}{sg}))$, 
and that the restriction of $\omega^n$ to $\cK(C_d(E^n))$ 
is a surjective map to $\cK(C_d(\s{E}{n}{sg}))$, 
whose kernel is $\cK(C_d(\s{E}{n}{rg}))$ (Lemma \ref{surjection4}). 
We denote by $\dot\pi_n^n\colon \F_T^{n}\to\cL(C_d(\s{E}{n}{sg}))$ 
the composition of the map $\pi_n^n\colon \F_T^{n}\to\cL(C_d(E^n))$ 
and the surjection $\omega^n\colon \cL(C_d(E^n))\to \cL(C_d(\s{E}{n}{sg}))$.

\begin{lemma}\label{Fn2}
For each $n\in\N$, 
we can define $*$-homomorphisms 
$\pi_k^{n}\colon \F_T^{n}\to\cL(C_d(\s{E}{k}{sg}))$ 
$(k=0,\ldots,n-1)$ 
such that $\bigcap_{k=0}^{n-1}\ker\pi_k^{n}=\G_T^n$. 
\end{lemma}

\begin{proof}
The proof goes by induction with respect to $n\in\N$. 
For $n=0$, we need to do nothing because $\G_T^0=\F_T^0$.
Assume that we have $*$-homomorphisms 
$\pi_k^{n}\colon \F_T^{n}\to\cL(C_d(\s{E}{k}{sg}))$ for $k=0,\ldots,n-1$ 
such that $\bigcap_{k=0}^{n-1}\ker\pi_k^{n}=\G_T^n$. 
Then we have
\begin{align*}
\bigcap_{k=0}^{n-1}\ker\pi_k^{n}\cap\ker\dot\pi_n^n
&=\G_T^n\cap\ker\dot\pi_n^n
=\varPhi^n\big(\cK(C_d(E^n))\cap\ker\omega^n\big)\\
&=\varPhi^n\big(\cK(C_d(\s{E}{n}{rg}))\big)
=\G_T^{n}\cap \G_T^{n+1}.
\end{align*}
Hence the maps $\pi_0^{n},\ldots,\pi_{n-1}^{n}$ and $\dot\pi_n^n$ 
factor through maps
$$\widetilde{\pi}_k^{n+1}\colon 
\F_T^{n}/(\G_T^n\cap\G_T^{n+1})\to\cL(C_d(\s{E}{k}{sg}))
\quad (k=0,\ldots,n),$$
and we see that $\bigcap_{k=0}^{n}\ker\widetilde{\pi}_k^{n+1}=0$.
For $k=0,\ldots,n$, 
we define 
$\pi_k^{n+1}\colon \F_T^{n+1}\to\cL(C_d(\s{E}{k}{sg}))$ 
by the composition of
the quotient map $\F_T^{n+1}\to \F_T^{n+1}/\G_T^{n+1}$,
the isomorphism $\F_T^{n+1}/\G_T^{n+1}\cong\F_T^{n}/(\G_T^n\cap\G_T^{n+1})$, 
and $\widetilde{\pi}_k^{n+1}$.
Then we have $\bigcap_{k=0}^{n}\ker\pi_k^{n+1}=\G_T^{n+1}$.
\end{proof}

For $k=0,\ldots,n-1$, 
the $*$-homomorphism $\pi_k^{n}\colon \F_T^{n}\to\cL(C_d(\s{E}{k}{sg}))$ 
defined in the proof of Lemma \ref{Fn2} 
is determined by
$$\pi_k^{n}(x)
=\left\{\begin{array}{ll}
\dot{\pi}_k^k(x)& \mbox{ if }x\in \F_T^{k},\\
0& \mbox{ if }x\in \G_T^{l} \mbox{ for }k<l\leq n.
\end{array}\right.$$

\begin{proposition}\label{Fn}
For $n\in\N$, the map 
$$\bigoplus_{k=0}^{n}\pi_k^n\colon  \F_T^n \to 
\bigoplus_{k=0}^{n-1}\cL(C_d(\s{E}{k}{sg}))\oplus\cL(C_d(E^n))$$
is injective.
\end{proposition}

\begin{proof}
This follows from 
$$\ker\bigg(\bigoplus_{k=0}^{n}\pi_k^n\bigg)=\bigcap_{k=0}^{n}\ker\pi_k^n
=\G_T^n\cap\ker\pi_n^n=0.$$
\end{proof}

\begin{definition}
Let $E$ be a topological graph. 
A path $e=(e_n,\ldots,e_1)\in E^n$ for $n\geq 1$ is called 
a {\em loop} if $r(e)=d(e)$, 
and the vertex $r(e)=d(e)$ is called the {\em base point} of the loop $e$. 
A loop $e=(e_n,\ldots,e_1)$ is 
said to be {\em without entrances} 
if $r^{-1}(r(e_k))=\{e_{k}\}$ for $k=1,\ldots,n$.
\end{definition}

It is easy to see that when $E$ is an ordinary dynamical system, 
every loops are without entrances and 
$v\in E^0$ is a base point of a loop if and only if 
it is a periodic point. 
We generalize the notion of topological freeness 
of homeomorphisms to topological correspondences. 
Recall that a homeomorphism on the space is called topologically free
if the set of periodic points has an empty interior 
(see \cite{AS}, \cite{T1}, \cite{ELQ}). 

\begin{definition}
A topological graph $E$ 
is said to be {\em topologically free} 
if the set of base points of loops without entrances
has an empty interior. 
\end{definition}

One can easily see that topological freeness coincides 
with {\em Condition L} when a graph $E$ is discrete (see \cite{KPR}). 
We will show that when a topological graph $E$ is topologically free, 
$\rho_T$ is an isomorphism if and only if $T^0$ is injective 
(Theorem \ref{CKUT}). 
To do so, we need the following notion and many lemmas.

\begin{definition}
Let $n$ be an integer with $n\geq 1$. 
A path $e=(e_n,\ldots,e_1)\in E^n$ is said to be 
{\em returning} if $e_1=e_k$ for some $k\in\{2,\ldots,n\}$. 
Otherwise $e$ is said to be {\em non-returning}. 
A non-empty set $U\subset E^n$ 
is said to be {\em non-returning} if 
$e_1\neq e_k'$ for every $k=2,\ldots, n$ and every
$(e_m,\ldots,e_1),(e_m',\ldots,e_1')\in U$. 
\end{definition}

\begin{lemma}\label{ckut1}
Let $V$ be an open subset of $E^0$, 
and $e=(e_n,\ldots,e_1)\in E^n$ $(n\geq 1)$ 
be a non-returning path with $r^n(e)\in V$. 
Then there exists a non-empty open set 
$U\subset (r^n)^{-1}(V)\subset E^n$ 
which is non-returning. 
\end{lemma}

\begin{proof}
Take open subsets $U_1,\ldots,U_n$ of $E^1$ such that 
$e_k\in U_k$ for $k=1,\ldots,n$, 
$U_k\cap U_1=\emptyset$ for $k=2,\ldots,n$, 
and $r^1(U_n)\subset V$. 
Then a non-empty open subset 
$U=(U_n\times\cdots\times U_1)\cap E^n$ of $E^n$ 
is non-returning. 
\end{proof}

\begin{lemma}\label{ckut4}
Suppose that an open set $U\subset E^m$ 
is non-returning. 
Take $\zeta\in C_c(U)\subset C_d(E^m)$ 
and $\xi\in C_d(E^n)$ for $1\leq n\leq m-1$. 
Then we have $T^m(\zeta)^*T^n(\xi)T^m(\zeta)=0$ 
for any Toeplitz $E$-pair $T$. 
\end{lemma}

\begin{proof}
By Lemma \ref{k=m-n}, 
we have $T^m(\zeta)^*T^n(\xi)T^m(\zeta)=T^n(\eta)$
where $\eta\in C_d(E^n)$ is defined by 
$$\eta(e_n,\ldots,e_1)
=\hspace*{-0.5cm}\sum_{\stackrel
{\mbox{$\scriptstyle (e_{n+m},\ldots,e_{n+1})\in E^m$}}
{\mbox{$\scriptstyle d(e_{n+1})=r(e_n)$}}}\hspace*{-0.5cm}
\overline{\zeta(e_{n+m},\ldots,e_{n+1})}
\xi(e_{n+m},\ldots,e_{m+1})\zeta(e_m,\ldots,e_1).$$
For each $(e_n,\ldots,e_1)\in E^n$, 
we see that $(e_{n+m},\ldots,e_{n+1})\in U$ 
implies $(e_m,\ldots,e_1)\notin U$ 
by the assumption on $U$. 
Hence we have $\eta=0$. 
Thus $T^m(\zeta)^*T^n(\xi)T^m(\zeta)=0$. 
\end{proof}

\begin{lemma}\label{ckut3}
Suppose that an open subset $V$ of $E^0$ 
satisfies that $(r^k)^{-1}(V)\neq\emptyset$ and $(r^{k+1})^{-1}(V)=\emptyset$ 
for some $k\in\N$.
Take $\zeta\in C_c((r^k)^{-1}(V))\subset C_d(E^k)$ 
and $\xi\in C_d(E^l)$ for $l>k$. 
Then we have $T^k(\zeta)^*T^l(\xi)=0$ 
for any Toeplitz $E$-pair $T$. 
\end{lemma}

\begin{proof}
By Lemma \ref{k=m-n}, 
we have $T^k(\zeta)^*T^l(\xi)=T^m(\eta)$ where $m=l-k>0$ and 
$\eta\in C_d(E^m)$ is determined by 
$$\eta(e)=\!\sum_{\stackrel
{\mbox{$\scriptstyle e'\in E^k$}}
{\mbox{$\scriptstyle d^k(e')=r^m(e)$}}}\!\!
\overline{\zeta(e')}\xi(e',e)\qquad (e\in E^m).$$
For $e\in E^m$, 
there exists no $e'\in (r^k)^{-1}(V)\subset E^k$ with $d^k(e')=r^m(e)$ 
by the assumption. 
Hence $\eta=0$. 
Thus $T^k(\zeta)^*T^l(\xi)=0$. 
\end{proof}

\begin{lemma}\label{ckut2}
Suppose that a topological graph 
$E=(E^0,E^1,d,r)$ is topologically free. 
For an open subset $V$ of $E^0$ and a positive integer $n$, 
either $(r^n)^{-1}(V)=\emptyset$ or 
there exists a non-returning path $e\in E^m$ with $m\geq n$ 
such that $r^m(e)\in V$. 
\end{lemma}

\begin{proof}
To the contrary, 
assume that an open subset $V$ of $E^0$ 
satisfies that $(r^n)^{-1}(V)\neq\emptyset$ and
that every path in $(r^m)^{-1}(V)$ is returning 
for every $m\geq n$. 
Take $e=(e_n,\ldots,e_1)\in (r^n)^{-1}(V)$ arbitrarily. 
Since $e$ is returning, 
there exists $k_0$ with $2\leq k_0\leq n$ such that $e_{k_0}=e_1$.
We will show that $r^{-1}(r(e_l))=\{e_l\}$ for $l=1,\ldots,k_0-1$. 
To derive a contradiction,
assume that there exist an integer $l$ with $1\leq l<k_0$
and $e\in E^1$ such that $r(e)=r(e_l)$ and $e\neq e_l$. 
Set $v=r(e)=r(e_l)$. 
Let $k_1$ be a maximal integer satisfying $r(e_{k_1})=v$. 
We have $d(e_{k_1+1})=v$ and $e_k\neq e$ for $k=k_1+1,\ldots,n$.
Since $(e_{k_0-1},\ldots,e_1)$ is a loop, 
we can find a loop $e'=(e_m',\ldots,e_1')$ 
such that $e_m'=e_l$ and $r(e_k')\neq v$ for $k=1,\ldots,m-1$.
Hence we have $e_k'\neq e$ for $k=1,\ldots,m$.
Set 
$$e''=(e_n,\ldots,e_{k_1+1},e',e',\ldots,e',e)\in E^{n'}$$
where $e'$'s are repeated so that $n'\geq n$.
Then $e''\in E^{n'}$ is non-returning path with $r(e'')\in V$.
This is a contradiction. 
Therefore, 
we have shown that for every $(e_n,\ldots,e_1)\in (r^n)^{-1}(V)$, 
there exists an integer $k_0$ with $2\leq k_0\leq n$ such that 
$(e_{k_0-1},\ldots,e_1)$ is a loop without entrances.
Thus each element in a non-empty open subset $d^n((r^n)^{-1}(V))$ of $E^0$ 
is a base point of a loop without entrances. 
This contradicts the fact that $E$ is topologically free.
\end{proof}

\begin{proposition}\label{ckut5}
Let $E=(E^0,E^1,d,r)$ be a topologically free 
topological graph, 
and $T=(T^0,T^1)$ be an injective Cuntz-Krieger $E$-pair.
Take $n_l,m_l\in\N$, 
$\xi_l\in C_d(E^{n_l})$ and $\eta_l\in C_d(E^{m_l})$ 
for $l=1,2,\ldots,L$. 
Define 
$$x=\sum_{l=1}^LT^{n_l}(\xi_l)T^{m_l}(\eta_l)^*,\quad 
x_0=\sum_{n_l=m_l}T^{n_l}(\xi_l)T^{m_l}(\eta_l)^*.$$ 
Then for arbitrary $\e>0$ 
there exist $a,b\in C^*(T)$ and $f\in C_0(E^0)$ 
such that $\|a\|,\|b\|\leq 1$, 
$\|f\|=\|x_0\|$, and $\|a^*xb-T^0(f)\|<\e$. 
\end{proposition}

\begin{proof}
Set $n=\max\{n_1,\ldots,n_L,\ m_1,\ldots,m_L\}$.
We have $x_0\in\F_T^n$.
Since 
$$\bigoplus_{k=0}^{n}\pi_k^n\colon \F_T^n\to 
\bigoplus_{k=0}^{n-1}\cL(C_d(\s{E}{k}{sg}))\oplus\cL(C_d(E^n))$$
is injective, 
there exists an integer $k$ with $0\leq k\leq n$ 
such that $\|x_0\|=\|\pi_k^n(x_0)\|$. 

\case{{\em Case 1}, $k\leq n-1$.}
There exist $\xi',\eta'\in C_d(\s{E}{k}{sg})$ 
with $\|\xi'\|=\|\eta'\|=1$ 
such that 
$$\big\|\ip{\xi'}{\pi_k^n(x_0)\eta'}\big\|
>\|\pi_k^n(x_0)\|-\e=\|x_0\|-\e.$$
By Lemma \ref{surjection}, 
we can find $\xi,\eta\in C_d(E^k)$ with $\|\xi\|=\|\eta\|=1$ 
such that $\xi|_{\s{E}{k}{sg}}=\xi'$ and $\eta|_{\s{E}{k}{sg}}=\eta'$.
For each $l=1,\ldots,L$, 
there exist $n_l',m_l'\in\N$ 
and $\xi_l'\in C_d(E^{n_l'}),\ \eta_l'\in C_d(E^{m_l'})$ 
such that 
$$T^{n_l'}(\xi_l')T^{m_l'}(\eta_l')^* 
=T^k(\xi)^*T^{n_l}(\xi_l)T^{m_l}(\eta_l)^*T^k(\eta).$$
We set 
$$y=\sum_{l=1}^LT^{n_l'}(\xi_l')T^{m_l'}(\eta_l')^*\quad 
\mbox{and}\quad 
T^0(g)=\sum_{n_l'=m_l'=0}T^{n_l'}(\xi_l')T^{m_l'}(\eta_l')^*.$$ 
Then we have $y=T^k(\xi)^*xT^k(\eta)$ and 
$g|_{\s{E}{0}{sg}}=\ip{\xi'}{\pi_k^n(x_0)\eta'}$.
Since $\|x_0\|-\e<\|g|_{\s{E}{0}{sg}}\|\leq\|x_0\|$, 
we can find $v\in \s{E}{0}{sg}$ with $\|x_0\|-\e<|g(v)|\leq\|x_0\|$.
Since $\s{E}{0}{sg}=\overline{\s{E}{0}{sce}}\cup\s{E}{0}{inf}$,
there are two cases, 
namely the case that $v\in \overline{\s{E}{0}{sce}}$ 
and the case that $v\in \s{E}{0}{inf}$.

\case{{\em Subcase 1.1}, $v\in\overline{\s{E}{0}{sce}}$.}
There exists $v'\in \s{E}{0}{sce}$ with $\|x_0\|-\e<|g(v')|\leq\|x_0\|$.
Since $v'\in \s{E}{0}{sce}$, 
we can find a neighborhood $V$ of $v'$ such that $r^{-1}(V)=\emptyset$
and $|g(v'')|<\|x_0\|+\e$ for all $v''\in V$.
Take $h\in C_0(V)$ with $0\leq h\leq 1$ and $h(v')=1$. 
We set $a=T^k(\xi)T^0(h)$ and $b=T^k(\eta)T^0(h)$. 
Then we have $\|a\|,\|b\|\leq 1$ 
and $a^*xb=T^0(h)yT^0(h)=T^0(h)T^0(g)T^0(h)$ 
because if either $n_l'$ or $m_l'$ is grater than $0$, 
then $T^0(h)T^{n_l'}(\xi_l')T^{m_l'}(\eta_l')^*T^0(h)=0$
by Lemma \ref{ckut3}. 
Define $f'=hgh$. 
Then we have $\|x_0\|-\e<\|f'\|<\|x_0\|+\e$. 
Set $f=\|x_0\|f'/\|f'\|$. 
Then we have $\|f\|=\|x_0\|$ and $\|f-f'\|=\big|\|x_0\|-\|f'\|\big|<\e$.
Thus we get 
$$\|a^*xb-T^0(f)\|=\|T^0(f')-T^0(f)\|<\e.$$

\case{{\em Subcase 1.2}, $v\in \s{E}{0}{inf}$.}
Choose a positive number $\delta$ such that $\delta<\e$ and 
$\|x_0\|-\delta<|g(v)|$. 
There exists a neighborhood $V$ of $v$ 
such that $\|x_0\|-\delta<|g(v')|<\|x_0\|+\delta$ for all $v'\in V$.
For $l$ with $n_l'=0$, we set $\xi_l''=\xi_l'$. 
For $l$ with $n_l'\geq 1$, we set $\xi_l''\in C_c(E^{n_l'})$ 
sufficiently close to $\xi_l'$. 
Similarly, we set $\eta_l''=\eta_l'$ for $l$ with $m_l'=0$, 
and set $\eta_l''\in C_c(E^{n_l'})$ 
sufficiently close to $\eta_l'$ for $l$ with $m_l'\geq 1$. 
We set $y'=\sum_{l=1}^LT^{n_l'}(\xi_l'')T^{m_l'}(\eta_l'')^*$.
We can choose $\xi_l''$ and $\eta_l''$ so that $\|y-y'\|<\e-\delta$. 
We have $T^0(g)=\sum_{n_l'=m_l'=0}T^{n_l'}(\xi_l'')T^{m_l'}(\eta_l'')^*$.
Set 
\begin{align*}
K=\bigcup_{n_l'\geq 1}&\big\{e_{n_l'}\in E^1\ \big| \mbox{ there exists }
e=(e_{n_l'},\ldots,e_1)\in\supp(\xi_l'')\big\}\\
&\cup\bigcup_{m_l'\geq 1}\big\{e_{m_l'}\in E^1\ \big| \mbox{ there exists }
e=(e_{m_l'},\ldots,e_1)\in\supp(\eta_l'')\big\},
\end{align*}
which is a compact subset of $E^1$.
Since $v\notin \s{E}{0}{fin}$, we have $r^{-1}(V)\setminus K\neq\emptyset$.
Take an open set $U\subset r^{-1}(V)\setminus K$ 
such that the restriction of $d$ to $U$ is injective.
Let $\zeta\in C_c(U)\subset C_d(E^1)$ be an element 
with $0\leq\zeta\leq 1$ and $\zeta(e)=1$ for some $e\in U$.
We have $\|\zeta\|=1$. 
We get $T^1(\zeta)^*T^{n_l'}(\xi_l'')T^{m_l'}(\eta_l'')^*T^1(\zeta)=0$
if either $n_l'$ or $m_l'$ is not zero. 
We set $a=T^k(\xi)T^1(\zeta)$, $b=T^k(\eta)T^1(\zeta)$ 
and $f'=\ip{\zeta}{\pi_r(g)\zeta}$. 
Then we have $\|a\|,\|b\|\leq 1$ and 
\begin{align*}
\|a^*xb-T^0(f')\|
&=\|T^1(\zeta)^*yT^1(\zeta)-T^1(\zeta)^*T^0(g)T^1(\zeta)\|\\
&=\|T^1(\zeta)^*yT^1(\zeta)-T^1(\zeta)^*y'T^1(\zeta)\|
\leq\|y-y'\|<\e-\delta.
\end{align*}
For $e'\in U$, 
we have 
$$|f'(d(e'))|=|\overline{\zeta(e')}g(r(e'))\zeta(e')|
\leq|g(r(e'))|<\|x_0\|+\delta$$
because the restriction of $d$ to $U$ is injective. 
We also have 
$$|f'(d(e))|=|\overline{\zeta(e)}g(r(e))\zeta(e)|=|g(r(e))|
>\|x_0\|-\delta.$$ 
Hence we get $\big|\|f'\|-\|x_0\|\big|<\delta$. 
Therefore $f=\|x_0\|f'/\|f'\|$ satisfies 
that $\|f\|=\|x_0\|$ and $\|f-f'\|<\delta$.
Thus we have 
$$\|a^*xb-T^0(f)\|\leq\|a^*xb-T^0(f')\|+\|f'-f\|<\e.$$

\case{{\em Case 2}, $k=n$.}
Next we consider the case that $\|x_0\|=\|\pi_n^n(x_0)\|$. 
We can find $\xi,\eta\in C_d(E^n)$ with $\|\xi\|=\|\eta\|=1$ 
such that $\|\ip{\xi}{\pi_n^n(x_0)\eta}\|>\|x_0\|-\e$. 
Set $g=\ip{\xi}{\pi_n^n(x_0)\eta}\in C_0(E^0)$.
There exists a non-empty open set $V$ of $E^0$ such that
$|g(v)|>\|x_0\|-\e$ for $v\in V$. 
When $n_l>m_l$, 
we have 
$$T^n(\xi)^*T^{n_l}(\xi_l)T^{m_l}(\eta_l)^*T^n(\eta)=T^{n_l'}(\xi_l')$$
for some $\xi_l'\in C_d(E^{n_l'})$ where $n_l'=n_l-m_l$.
Similarly when $n_l<m_l$,
we have 
$$T^n(\xi)^*T^{n_l}(\xi_l)T^{m_l}(\eta_l)^*T^n(\eta)=T^{m_l'}(\eta_l')^*$$
for some $\eta_l'\in C_d(E^{m_l'})$ where $m_l'=m_l-n_l$.
We have 
$$T^n(\xi)^*xT^n(\eta)=T^0(g)+\sum_{n_l>m_l}T^{n_l'}(\xi_l')
+\sum_{n_l<m_l}T^{m_l'}(\eta_l')^*.$$

\case{{\em Subcase 2.1}, $(r^{n+1})^{-1}(V)=\emptyset$.}
Take an integer $k$ with $0\leq k\leq n$ 
satisfying that $(r^k)^{-1}(V)\neq\emptyset$
and $(r^{k+1})^{-1}(V)=\emptyset$.
Take $\zeta\in C_c((r^k)^{-1}(V))\subset C_d(E^k)$ such that 
$\|\zeta\|=1$ and $\zeta(e)=1$ for some $e\in (r^k)^{-1}(V)$.
Set $a=T^n(\xi)T^k(\zeta)$ and $b=T^n(\eta)T^k(\zeta)$. 
Then we have $\|a\|,\|b\|\leq 1$. 
We see that 
$T^k(\zeta)^*T^{n_l'}(\xi_l')T^k(\zeta)=0$ and 
$T^k(\zeta)^*T^{m_l'}(\eta_l')^*T^k(\zeta)=0$
by Lemma \ref{ckut3}.
Hence we get 
$$a^*xb=T^k(\zeta)^*T^0(g)T^k(\zeta).$$
Set $f'=\ip{\zeta}{\pi_{r^k}(g)\zeta}$.
Then we have $\|f'\|\leq\|g\|\leq\|x_0\|$ and
$$\|f'\|\geq |f'(d(e))|
\geq |\overline{\zeta(e)}g(r^k(e))\zeta(e)|\geq \|x_0\|-\e.$$
Therefore we get $\big|\|f'\|-\|x_0\|\big|<\e$. 
Hence $f=\|x_0\|f'/\|f'\|$ satisfies 
that $\|f\|=\|x_0\|$ and $\|a^*xb-T^0(f)\|<\e$. 

\case{{\em Subcase 2.2}, $(r^{n+1})^{-1}(V)\neq\emptyset$.}
By Lemma \ref{ckut2}, 
there exists a non-returning path $e\in E^m$ 
with $m\geq n+1$ and $r(e)\in V$. 
By Lemma \ref{ckut1}, 
we can find a non-empty open set $U\subset (r^m)^{-1}(V)\subset E^m$
which is non-returning. 
Choose $\zeta\in C_c(U)\subset C_d(E^m)$ such that 
$\|\zeta\|=1$ and $\zeta(e)=1$ for some $e\in U$. 
Set $a=T^n(\xi)T^m(\zeta)$, $b=T^n(\eta)T^m(\zeta)$ 
and $f'=\ip{\zeta}{\pi_{r^m}(g)\zeta}$. 
Then we have $\|a\|,\|b\|\leq 1$ and 
$a^*xb=T^m(\zeta)^*yT^m(\zeta)=T^0(f')$ 
by Lemma \ref{ckut4}.
Similarly as the proof in Subcase 2.1, 
we have $\big|\|f'\|-\|x_0\|\big|<\e$. 
Hence $f=\|x_0\|f'/\|f'\|$ satisfies 
that $\|f\|=\|x_0\|$ and $\|a^*xb-T^0(f)\|<\e$. 

The proof is completed.
\end{proof}

\begin{remark}
For a positive element $x\in \cL(X)$ where $X$ is a Hilbert module, 
we have $\ip{\xi}{x\xi}\geq 0$ for any $\xi\in X$ and 
$\|x\|=\sup_{\|\xi\|=1}\|\ip{\xi}{x\xi}\|$. 
Therefore 
if we further assume $x_0\geq 0$ 
in the assumption of Proposition \ref{ckut5},
then we can take $a,b\in C^*(T)$ and $f\in C_0(E^0)$ 
in the conclusion there so that $a=b$ and $f\geq 0$
because we can take $\xi=\eta$ 
in both Case 1 and Case 2 in the above proof. 
We will use this fact in \cite{Ka4}. 
\end{remark}

The following is our version of the Cuntz-Krieger uniqueness theorem.

\begin{theorem}\label{CKUT}
If a topological graph $E=(E^0,E^1,d,r)$ is topologically free, 
then the map $\rho_{T}\colon \cO(E)\to C^*(T)$ is an isomorphism
for any injective Cuntz-Krieger $E$-pair $T=(T^0,T^1)$.
\end{theorem}

\begin{proof}
We check the condition (iii) in Theorem \ref{GIUT}. 
By Proposition \ref{ckut5}, 
for any $\e>0$ we have $\|x_0\|=\|f\|\leq\|a^*xb\|+\e\leq\|x\|+\e$ 
for $x=\sum_{l=1}^LT^{n_l}(\xi_l)T^{m_l}(\eta_l)^*$ 
and $x_0=\sum_{n_l=m_l}T^{n_l}(\xi_l)T^{m_l}(\eta_l)^*$ 
where $\xi_l\in C_d(E^{n_l})$ and $\eta_l\in C_d(E^{m_l})$. 
Hence $x\mapsto x_0$ gives us a well-defined norm-decreasing 
linear map from 
$$\spa\{T^n(\xi)T^m(\eta)^*\mid 
\xi\in C_d(E^n),\ \eta\in C_d(E^m),\ n,m\in\N\},$$
to $\F_T\subset C^*(T)$. 
It extends a linear map $\varPsi_T$ from $C^*(T)$ to $\F_T$ 
which is identity on $\F_T$. 
Hence $\varPsi_T$ is a conditional expectation onto $\F_T$ such that 
$\varPsi_T(T^n(\xi)T^m(\eta)^*)=\delta_{n,m}T^n(\xi)T^m(\eta)^*$ 
for $\xi\in C_d(E^n),\ \eta\in C_d(E^m)$. 
By Theorem \ref{GIUT}, $\rho_{T}$ is an isomorphism. 
\end{proof}

\begin{remark}
In the theorem above, 
the assumption that a topological graph $E$ is 
topologically free is needed. 
When $E$ is not topologically free, 
there exists an injective Cuntz-Krieger $E$-pair $T$ 
such that the map $\rho_{T}\colon \cO(E)\to C^*(T)$ is not an isomorphism 
(see \cite{Ka3}). 
\end{remark}

\section{$KK$-groups of $\cO(E)$}\label{Kgroups}

In this section, 
we prove that our $C^*$-algebras $\cO(E)$ are always nuclear 
and satisfy the Universal Coefficient Theorem (UCT) of \cite{RoSc}, 
and compute their $KK$-groups. 
To this end, we need a short exact sequence 
which is (almost) established in Section \ref{Fock}. 
Take a topological graph $E=(E^0,E^1,d,r)$. 
In Section \ref{Fock}, 
we defined a Toeplitz $E$-pair $\sigma$ 
on $\cL(C_d(E^{*}))$ 
which was shown to be universal (Proposition \ref{univToep}). 
Hence we can identify $C^*(\sigma)$ with $\cT(E)$. 
The $C^*$-algebra $\cT(E)=C^*(\sigma)$ has an ideal $\cK(C_d(\s{E}{*}{rg}))$ 
which we will denote $I_E$. 
The quotient $\cT(E)/I_E$ 
is the $C^*$-algebra $C^*(\tau)$ 
which was shown to be isomorphic to 
the $C^*$-algebra $\cO(E)$ in Corollary \ref{Fockrep}. 
Hence we have a following short exact sequence:
$$\begin{CD}
0 @>>> I_E@>j>> \cT(E)@>>> \cO(E)@>>> 0.
\end{CD}$$

\begin{proposition}
The $C^*$-algebra $\cO(E)$ is nuclear. 
\end{proposition}

\begin{proof}
The $C^*$-algebra $\cT(E)$ 
is isomorphic to the augmented Toeplitz algebra $\widetilde{\cT}_{C_d(E^1)}$.
There is 
a folklore\footnote{In \cite{Ka5}, we give a proof of this folklore.} that 
an augmented Toeplitz algebra $\widetilde{\cT}_{X}$ 
of a $C^*$-correspondences $X$ over $A$ is nuclear 
if and only if $A$ is nuclear 
(the proof goes similarly as \cite{DS}). 
Hence the $C^*$-algebra $\widetilde{\cT}_{C_d(E^1)}$ is nuclear 
because $C_0(E^0)$ is nuclear. 
Since nuclearity inherits to quotients 
(see, for example, \cite[Corollary 2.5]{W}), 
the $C^*$-algebra $\cO(E)$ is nuclear. 
\end{proof}

It is well-known that a $C^*$-algebra $C_0(E^0)$ is separable 
if and only if $E^0$ is second countable. 
Similarly we have the following. 

\begin{lemma}\label{sep2ndctbl}
The Banach space $C_d(E^1)$ is separable 
if and only if $E^1$ is second countable. 
\end{lemma}

\begin{proof}
If $\{\xi_k\}_{k\in\N}$ is a countable dense set of $C_d(E^1)$, 
then $\{U_k\}_{k\in\N}$ is a countable open basis of $E^1$ 
where $U_k=\{e\in E^1\mid |\xi_k(e)|<1\}$. 
Hence if $C_d(E^1)$ is separable, 
$E^1$ is second countable. 
Conversely if $E^1$ is second countable, 
we can find a countable subset $X$ of $C_c(E^1)$ 
such that for every open subset $U$ of $E^1$, 
every elements of $C_c(U)$ 
can be uniformly approximated by elements of $X\cap C_c(U)$. 
By Lemma \ref{dense2}, 
we see that the countable subset $X$ is dense in $C_d(E^1)$ 
with respect to the norm $\|\cdot\|$. 
Hence $C_d(E^1)$ is separable. 
\end{proof}

We say that a topological graph 
$E=(E^0,E^1,d,r)$ is {\em second countable} 
if both $E^0$ and $E^1$ are second countable. 

\begin{proposition}
The $C^*$-algebra $\cO(E)$ is separable 
if and only if $E$ is second countable. 
\end{proposition}

\begin{proof}
If $\cO(E)$ is separable, 
then both $C_0(E^0)$ and $C_d(E^1)$ are separable 
because $t^0$ and $t^1$ are isometric. 
Conversely if both $C_0(E^0)$ and $C_d(E^1)$ are separable 
then $\cO(E)$ is separable 
because $\cO(E)$ is generated by the images of $C_0(E^0)$ and $C_d(E^1)$. 
Now the proof ends by Lemma \ref{sep2ndctbl}. 
\end{proof}

In the rest of this section, 
we assume that $E$ is second countable. 

\begin{lemma}
The Hilbert $C_0(\s{E}{0}{rg})$-module $C_d(\s{E}{*}{rg})$ is full. 
Hence it gives a strong Morita equivalence 
between $I_E=\cK(C_d(\s{E}{*}{rg}))$
and $C_0(\s{E}{0}{rg})$. 
\end{lemma}

\begin{proof}
It is easy to verify. 
\end{proof}

We denote by $[C_d(\s{E}{*}{rg})]\in KK(I_E,C_0(\s{E}{0}{rg}))$ 
the element defined by the imprimitivity bimodule $C_d(\s{E}{*}{rg})$. 
This element gives a $KK$-equivalence of $I_E$ and $C_0(\s{E}{0}{rg})$. 
Next we see that the inclusion $\sigma^0\colon C_0(E^0)\to \cT(E)$ gives 
a $KK$-equivalence from $C_0(E^0)$ to $\cT(E)$
following \cite{Pi}. 
We define a graded Kasparov module 
$(C_d(E^*)\oplus C_d(E^*), \pi\oplus\pi_+, T)$ as follows. 
A Toeplitz $E$-pair $\sigma$ on $\cL(C_d(E^*))$
gives an injective map $\pi\colon \cT(E)\to\cL(C_d(E^*))$. 
Define $E_+^*=\coprod_{n=1}^\infty E^n$ 
which is an open and closed subset of $E^*$. 
Hence $C_d(E_+^*)$ is a direct summand of $C_d(E^*)$ 
and we can consider $\cL(C_d(E_+^*))\subset \cL(C_d(E^*))$. 
This summand $C_d(E_+^*)$ is closed under two maps $\sigma^0,\sigma^1$. 
We write the restrictions of $\sigma^0,\sigma^1$ to $C_d(E_+^*)$ 
by $\sigma_+^0,\sigma_+^1$ respectively. 
Since the pair $\sigma_+=(\sigma_+^0,\sigma_+^1)$ 
is a restriction of 
a Toeplitz $E$-pair $\sigma=(\sigma^0,\sigma^1)$, 
it is also a Toeplitz $E$-pair on $\cL(C_d(E_+^*))$. 
This gives a $*$-homomorphism 
$\pi_+\colon \cT(E)\to\cL(C_d(E_+^*))$. 
Since $\cL(C_d(E_+^*))\subset \cL(C_d(E^*))$, 
we regard $\pi_+$ as a $*$-homomorphism to $\cL(C_d(E^*))$. 
Define an odd operator $T\in\cL(C_d(E^*)\oplus C_d(E^*))$ 
by $T(\xi\oplus \eta)=\eta\oplus \xi$. 
By \cite[Lemma 4.2]{Pi}, the triple 
$(C_d(E^*)\oplus C_d(E^*), \pi\oplus\pi_+, T)$ 
is a Kasparov module, 
and so it gives an element $\beta$ in 
$KK(\cT(E),C_0(E^0))$. 

\begin{lemma}
The element $\beta\in KK(\cT(E),C_0(E^0))$ 
is the inverse of 
the element $\sigma^0_*\in KK(C_0(E^0),\cT(E))$. 
Hence $\cT(E)$ is $KK$-equivalent to $C_0(E^0)$. 
\end{lemma}

\begin{proof}
See \cite[Theorem 4.4]{Pi}. 
\end{proof}

\begin{proposition}
For a second countable topological graph $E$,
the $C^*$-algebra $\cO(E)$ satisfies the UCT. 
\end{proposition}

\begin{proof}
Since $I_E$ and $\cT(E)$ are $KK$-equivalent 
to commutative $C^*$-algebras $C_0(\s{E}{0}{rg})$ and $C_0(E^0)$ 
respectively, 
they satisfy the UCT.
Now ``two among three principle'' shows that 
$\cO(E)$ satisfies the UCT. 
\end{proof}

Finally we compute the $KK$-groups of $\cO(E)$ 
in terms of the topological graph $E$. 
To do so, we examine the element $j_*\in KK(I_E,\cT(E))$ 
defined by the inclusion $j\colon I_E\to \cT(E)$. 
We denote by $[\pi_r]$ 
the element of $KK(C_0(\s{E}{0}{rg}), C_0(E^0))$ 
defined by a triple $(C_d(E^1),\pi_r,0)$ 
and denote by $\iota\colon C_0(\s{E}{0}{rg})\to C_0(E^0)$ a natural embedding. 
We have the following. 

\begin{lemma}\label{maponKK}
With the notation above, 
we have 
$$[C_d(\s{E}{*}{rg})]\otimes_{C_0(\s{E}{0}{rg})}
(\iota_*-[\pi_r])=j_*\otimes_{\cT(E)}\beta$$
in $KK(I_E,C_0(E^0))$. 
\end{lemma}

\begin{proof}
The proof is exactly the same as in \cite[Lemma 4.7]{Pi}, 
hence we omit it. 
\end{proof}

\begin{remark}
When $E^1=E^0$ and $d=\id$, 
there exists a natural isomorphism 
$\cK(C_d(E^1))\cong C_0(E^0)$. 
Under this isomorphism, 
the restriction of $\pi_r$ to $C_0(\s{E}{0}{rg})$ 
coincides with a $*$-homomorphism 
$\hat{r}\colon C_0(\s{E}{0}{rg})\ni f\mapsto f\circ r\in C_0(E^0)$. 
We see that the element $[\pi_r]\in KK\big(C_0(\s{E}{0}{rg}), C_0(E^0)\big)$ 
is defined by this $*$-homomorphism $\hat{r}$. 
For a general topological graph $E$, 
$\cK(C_d(E^1))$ is strongly Morita equivalent to 
the ideal $C_0(d(E^1))$ of $C_0(E^0)$ 
via the Hilbert module $C_d(E^1)$. 
Thus we have an element of 
$[C_d(E^1)]\in KK\big(\cK(C_d(E^1)), C_0(E^0)\big)$. 
The element $[\pi_r]\in KK\big(C_0(\s{E}{0}{rg}), C_0(E^0)\big)$ is 
the Kasparov product of 
the element $(\pi_r)_*\in KK\big(C_0(\s{E}{0}{rg}), \cK(C_d(E^1))\big)$ 
induced by $\pi_r\colon C_0(\s{E}{0}{rg})\to\cK(C_d(E^1))$ 
and $[C_d(E^1)]$. 
\end{remark}

\begin{proposition}\label{6term}
Let $E$ be a second countable topological graph. 
For any separable $C^*$-algebra $B$, 
we have the following two exact sequences:
$$\begin{CD}
KK_0(B,C_0(\s{E}{0}{rg})) @>>\iota_*-[\pi_r]> KK_0(B,C_0(E^0)) 
@>>t^0_*>  KK_0(B,\cO(E)) \\
@AAA @. @VVV \\
KK_1(B,\cO(E)) @<t^0_*<< KK_1(B,C_0(E^0)) 
@<\iota_*-[\pi_r]<< KK_1(B,C_0(\s{E}{0}{rg}))
\end{CD}$$
and
$$\begin{CD}
KK_0(C_0(\s{E}{0}{rg}),B) @<<\iota_*-[\pi_r]< KK_0(C_0(E^0),B) 
@<<t^0_*< KK_0(\cO(E),B) \\
@VVV @. @AAA \\
KK_1(\cO(E),B) @>t^0_*>> KK_1(C_0(E^0),B) 
@>\iota_*-[\pi_r]>> KK_1(C_0(\s{E}{0}{rg}),B).
\end{CD}$$
\end{proposition}

\begin{proof}
From the short exact sequence 
$$\begin{CD}
0 @>>> I_E@>j>> \cT(E)@>>> \cO(E)@>>> 0,
\end{CD}$$
we have two 6-term exact sequences of $KK$-groups
because $\cO(E)$ is nuclear. 
As we saw above, $I_E$ and $\cT(E)$ are $KK$-equivalent 
to $C_0(\s{E}{0}{rg})$ and $C_0(E^0)$ respectively. 
Under these isomorphisms, 
the element $j_*\in KK(I_E,\cT(E))$ 
coincides with $\iota_*-[\pi_r]\in KK(C_0(\s{E}{0}{rg}),C_0(E^0))$ 
by Lemma \ref{maponKK}. 
The composition of the map $\sigma^0\colon C_0(E^0)\to\cT(E)$ 
and the quotient map $\cT(E)\to\cO(E)$ is $t^0\colon C_0(E^0)\to\cO(E)$. 
Now we have the desired 6-term exact sequences. 
\end{proof}

\begin{corollary}\label{Kgroup}
For a second countable\footnote{In \cite{Ka5}, 
we prove this corollary without the assumption of second countability.} 
topological graph $E$, 
we have the following exact sequence of $K$-groups: 
$$\begin{CD}
K_0(C_0(\s{E}{0}{rg})) @>>\iota_*-[\pi_r]> K_0(C_0(E^0)) 
@>>t^0_*>  K_0(\cO(E)) \\
@AAA @. @VVV \\
K_1(\cO(E)) @<t^0_*<< K_1(C_0(E^0)) 
@<\iota_*-[\pi_r]<< K_1(C_0(\s{E}{0}{rg})).
\end{CD}$$
\end{corollary}

Finally we give a new proof of the computation of $K$-groups 
of graph algebras by using Corollary \ref{Kgroup}. 
Let $E=(E^0,E^1,d,r)$ be a discrete graph. 
The group $K_0(C_0(E^0))$ is isomorphic to
a free abelian group $\Z^{E^0}$ 
whose generators are $\{[\delta_v]\}_{v\in E^0}$, 
and $K_1(C_0(E^0))=0$. 
We also have $K_0(C_0(\s{E}{0}{rg}))\cong\Z^{\s{E}{0}{rg}}$ 
and $K_1(C_0(\s{E}{0}{rg}))=0$. 
For $v\in \s{E}{0}{rg}$, 
we have $\iota_*([\delta_v])=[\delta_v]$. 
We will compute $[\pi_r]([\delta_v])\in K_0(C_0(E^0))$. 
By the computation done in Section \ref{graphalg}, 
we have $\pi_r(\delta_v)=\sum_{e\in r^{-1}(v)}\theta_{\delta_e,\delta_e}$ 
for $v\in \s{E}{0}{rg}$ 
(note that $r^{-1}(v)$ is a non-empty finite set). 
There exists an isomorphism 
$$\psi\colon K_0\big(\cK(C_d(E^1))\big)\to 
K_0\big(C_0(d(E^1))\big)\cong\Z^{d(E^1)}\subset \Z^{E^0}$$
defined by $\psi([\theta_{\delta_e,\delta_e}])=[\delta_{d(e)}]$
for $e\in E^1$. 
Using this map, we have 
$$[\pi_r]([\delta_v])=\psi([\pi_r(\delta_v)])
=\psi\bigg(\sum_{e\in r^{-1}(v)}[\theta_{\delta_e,\delta_e}]\bigg)
=\sum_{e\in r^{-1}(v)}[\delta_{d(e)}]\in K_0(C_0(E^0)).$$

Now Corollary \ref{Kgroup} gives us the following.

\begin{proposition}[{\cite[Proposition 2]{Sz}}, {\cite[Theorem 3.1]{DT2}}]
Let $E=(E^0,E^1,d,r)$ be a discrete graph. 
Define a homomorphism $\Delta\colon \Z^{\s{E}{0}{rg}}\to\Z^{E^0}$ 
by $\Delta([\delta_v])=[\delta_v]-\sum_{e\in r^{-1}(v)}[\delta_{d(e)}]$. 
Then we have isomorphisms $K_0(\cO(E))\cong \coker\Delta$ 
and $K_1(\cO(E))\cong\ker\Delta$. 
\end{proposition}

\end{document}